
%

%
\input amstex
\documentstyle{amsppt}
\loadbold
\def\cstar{$C^*$-algebra}
\def\esg{$E_0$-semigroup}

\def\<{\left<}										
\def\>{\right>}

\magnification=\magstep 1

\topmatter
\title Noncommutative flows I:\\
dynamical invariants
\endtitle

\author William Arveson
\endauthor

\affil Department of Mathematics\\
University of California\\Berkeley CA 94720, USA
\endaffil

\date 18 November 1995
\enddate
\thanks This research was supported by
NSF grants DMS92-43893 and DMS95-00291
\endthanks
\abstract
We show that a noncommutative dynamical system of the type
that occurs in quantum theory can often be associated with a
dynamical principle; that is, an infinitesimal structure
that completely determines the dynamics.
The nature of these dynamical principles
is similar to that of the second order differential
equations of classical mechanics, in that one can locate a space of
momentum operators, a ``Riemannian metric", and a potential.
These structures are classified in terms of geometric objects
which, in the simplest cases, occur in finite dimensional matrix
algebras.  As a consequence,
we obtain a new classification of \esg s acting on type $I$ factors.
\endabstract


\toc
\specialhead{Introduction}
\endspecialhead

\subhead 1.  Completely positive semigroups
\endsubhead
\subhead 2.  Differential operators on matrix algebras
\endsubhead
\subhead 3.  Elliptic operators
\endsubhead
\subhead 4.  Momentum
\endsubhead
\subhead 5.  Modular Cohomology
\endsubhead
\subhead 6.  Exactness
\endsubhead
\subhead 7.  Classification of elliptic operators
\endsubhead
\subhead 8.  Applications
\endsubhead
\specialhead{References}
\endspecialhead
\endtoc

\endtopmatter
\vfill\eject

\document

\subheading{Introduction}
The dynamical groups of quantum theory tend to be
one-parameter groups of automorphisms of \cstar s
generated by bounded functions of
certain observables.
Frequently, there is a distinguished state on the algebra
that is invariant under the action of the automorphism group
and which represents the vacuum.  Using a standard GNS construction,
one can arrange that the operator algebra acts on a Hilbert
space, that the automorphism group is implemented by a
one-parameter group of unitary operators on this representation
space, and that there is a unit vector in the Hilbert space,
fixed under the action of the unitary group, which gives
rise to the vacuum state.

We have settled on the following general context for confronting
dynamical systems of this type.

\proclaim{Definition A}
A {\rm history} is a triple $(M_{-}, U, v)$ consisting of a
one-paramter group of unitary operators $U = \{U_t: t\in \Bbb R\}$
acting on a separable Hilbert space $H$, a unit vector $v\in H$
satisfying
$$
U_t v = v,\qquad t\in \Bbb R,
$$
and a von Neumann algebra $M_{-}\subseteq \Cal B(H)$ such that
$$
\gamma_t(M_{-})\subseteq M_{-},\qquad t\leq 0,
$$
$\gamma = \{\gamma_t: t\in \Bbb R\}$ denoting the automorphism
group of $\Cal B(H)$ associated with $U$.
\endproclaim
The group $\gamma$ acts on $\Cal B(H)$ by
$$
\gamma_t(a) = U_t a U_t^*, \qquad t\in \Bbb R. \tag{0.1}
$$

There are two kinds of degeneracies that we
have found it necessary
to rule out.
Suppose first that $M_{-}$ is an abelian von Neumann algebra.
Then so is the larger von Neumann algebra
$$
M = (\bigcup_{t\geq 0} \gamma_t(M_{-}))^{\prime\prime}
$$
and $M$ is invariant under the full dynamical group
$\{\gamma_t: t\in \Bbb R\}$.  After introducing appropriate
coordinates we may realize $M$ as the algebra of all bounded random
variables on a standard probability space $(\Omega, P)$, in such
a way that the group action $\gamma$ is represented by a
one-parameter group of measure-preserving transformations
of $(\Omega, P)$.  $M_{-}$ becomes a subalgebra of $L^\infty(\Omega, P)$
that is invariant under negative time translations.

Now assume further that there is a von Neumann algebra
$M_0\subseteq M_{-}$ such that $M_{-}$ is generated by the
negative time translations of $M_0$:
$$
M_{-} = (\bigcup_{t\leq 0} \gamma_t(M_0))^{\prime\prime}.\tag{0.2}
$$
In this case we can find a single random variable
$$
X_0:\Omega \to \Bbb R
$$
with the property that $M_0$ consists of all bounded
measurable functions of $X_0$ (we may choose $X_0$
to be bounded if we wish).  We can define a stationary
stochastic process $\{X_t: t\in \Bbb R\}$ in the obvious
way
$$
X_t = \gamma_t(X_0), \qquad t\in \Bbb R,
$$
and thus $M_{-}$ appears as the ``past" of the
process $\{X_t: t\in \Bbb R\}$.  Thus this case falls completely
into the domain of stationary stochastic processes.

In general it may be impossible to find a ``sharp time"
subalgebra $M_0$ which satisfies (0.2).  Nevertheless, that
situation too is familiar in modern probability theory.
For example, that would be the the case when $M_{-}$ is
the ``past" subalgebra associated with a stationary random
distribution such as white noise.

For our purposes this probabilistic case is a degenerate one,
and we rule it out by imposing the following requirement
that forces us to the opposite extreme:

\proclaim{Definition B}
A History $(M_{-}, U, v)$ is called {\rm primary} if
$M_{-}$ is a factor.
\endproclaim

The second type of degeneracy represents a kind of {\it determinism}.
Let  $(M_{-}, U, v)$ be a history.  Then $M_{-}$ is invariant
under $\gamma_t$ for every $t\leq 0$, and if there is a single
$t_0<0$ for which $\gamma_{t_0}(M_{-}) = M_{-}$, then
$\gamma_t(M_{-}) = M_{-}$ for every negative $t$ and in fact
$M_{-}$ is invariant under the full dynamical group
$\{\gamma_t: t\in \Bbb R\}$.  One may interpret this as a
deterministic situation in which the future is a function
of the past (more precisely, every observable is a Borel
function of an observable belonging to $M_{-}$).
{}From our point of view, however, this
simply represents a bad choice of the subalgebra $M_{-}$.

These considerations have led us to
the following program for understanding
the nature of noncommutative dynamical systems.  One should seek
to classify primary histories $(M_{-}, U, v)$ which satisfy the
condition
$$
\gamma_t(M_{-})\subsetneqq M_{-}, \qquad \text{ for every } t<0.
$$
There are two significant equivalence relations for histories.
The first is an obvious notion of conjugacy: $(M_{-},U,v)$
and $(\tilde M_{-},\tilde U,\tilde v)$ are said to be
{\it conjugate} if there is a unitary operator $W: H\to \tilde H$
such that
$$
\align
Wv &= \tilde v \tag{i}\\
WU_t &= \tilde U_tW, \qquad t\in R \tag{ii}\\
WM_{-}W^{-1} &= \tilde M_{-}.  \tag{iii}
\endalign
$$
The second equivalence relation is {\it cocycle conjugacy},
and will be discussed in a subsequent paper of this series.

There are several ways of finding histories among other objects
that arise naturally in quantum physics.  For example, suppose that
we start quite simply
with a pair $(U, v)$ consisting of a one-parameter unitary
group $U$ acting on a Hilbert space $H$ and a unit vector
$v\in H$ such that $U_tv = v$ for every $t\in \Bbb R$.  Assume
further that $K$ is another Hilbert space and we are given
a representation $W_0$ of the canonical commutation relations
over $K$.  Thus
$$
f\in K \mapsto W_0(f) \in \Cal B(H)
$$
is a strongly continuous mapping from $K$ to the unitary
operators on $H$ for which
$$
W_0(f)W_0(g) = e^{i\omega(f,g)}W_0(f+g),\qquad f,g\in K,
$$
$\omega(f,g)$ denoting the imaginary part of the inner product
$\<f,g\>$ in $K$.  Then for every $t\in \Bbb R$ we can define
a representation $W_t: K\to \Cal B(H)$ of the CCRs by
$W_t(f) = U_tW_0(f) U_t^*$.  One may
think of the family $\{W_t: t\in \Bbb R\}$ as a noncommutative
quantum process.  In this case we have a natural candidate
for $M_{-}$:
$$
M_{-} = (\bigcup_{t\leq 0}W_t(K))^{\prime\prime},
$$
and thus we have a history $(M_{-}, U, v)$.

In this example there is a ``sharp time" von Neumann algebra associated
with time zero
$$
M_0 = W_0(K)^{\prime\prime}
$$
with the property that $M_{-}$ is generated by $M_0$ and its
negative time images under the group of automorphisms
associated with $U$.  In general, one cannot expect to have sharp
time algebras as this example does.  Nevertheless, in very
general circumstances (e.g., in any case where one has a quantum
system obeying the Haag-Kastler axioms \cite{16, pp. 106--107}) one can
find significant examples of histories $(M_{-}, U, v)$.

This paper deals with the classification of histories
$(M_{-},U,v)$ of the simplest type, namely
$$
\align
M_{-} \text{ is a factor of type $I$, and} \tag{0.3.1}\\
\gamma_t(M_{-})\subsetneqq M_{-}, \text{ for every }t<0.  \tag{0.3.2}
\endalign
$$
This leads to the problem of classifying {\it pairs} of
\esg s in the following way.

Let $M_{+}$ be the commutant of
$M_{-}$.  Then we obtain a pair of \esg s $\alpha^+$, $\alpha^-$
acting respectively on $M_+$ and $M_{-}$ by
$$
\align
\alpha^+_t(a) &= U_taU_t^*,\qquad t\geq 0, a\in M_+\\
\alpha^-_t(b) &= U_t^*bU_t,\qquad t\geq 0, b\in M_-.
\endalign
$$
Moreover, the vector state of $\Cal B(H)$ defined by $v$,
$\omega(a) = \<av,v\>$, can be restricted to
these two subalgebras to give normal states $\rho^+$, $\rho^-$
which are invariant under the actions of
the respective \esg s.  Thus the classification problem for
type $I$ histories reduces to the problem of classifying
{\it pairs} of \esg s which have distinguished normal invariant
states.  The theory of \esg s was initiated by Powers in \cite{24},
and has undergone considerable development during the past several years, see
\cite{1}, \cite{3}, \cite{4}, \cite{5}, \cite{25}, \cite{26}, \cite{27}
and references cited in these papers.

Notice too that while $\omega(a) = \<av,v\>$ is a pure state of
$\Cal B(H)$ its restrictions to the subfactors $M_+$ and
$M_-$ are {\it not} necessarily pure states on these algebras.
Thus we are led
to the problem of classifying pairs $(\alpha,\rho)$ where
$\alpha = \{\alpha_t: t\geq 0\}$ is an \esg\ acting on a
a type $I_\infty$ factor (which we may take as $\Cal B(H_\alpha)$
for some separable Hilbert space $H_\alpha$) and $\rho$ is
a normal state on $\Cal B(H_\alpha)$ satisfying
$$
\rho(\alpha_t(a)) = \rho(a), \qquad a\in \Cal B(H_\alpha), t\geq 0.
$$
Two such pairs $(\alpha,\rho)$ and $(\tilde \alpha,\tilde \rho)$
are said to be {\it conjugate} if there exists a $*$-isomorphism
$\theta : \Cal B(H_\alpha)\to \Cal B(H_{\tilde \alpha})$ such that
$$
\align
\tilde \alpha_t \circ \theta &= \theta \circ \alpha_t,\qquad t\geq 0,\\
\tilde \rho \circ \theta &= \rho.
\endalign
$$

Our main results below apply to the case where where $\rho$ is {\it not}
a pure state and is weakly continuous.  Equivalently, if we realize
$\rho$ in the form
$$
\rho(a) = \text{trace}(\Omega a)
$$
where $\Omega$ is a positive trace-class operator on $H_\alpha$,
then we are assuming  that
$$
2\leq \text{ rank } \Omega < \infty.  \tag{0.4}
$$
In the extreme case where $\Omega$ is of rank $1$, $\rho$ is
a pure state of $\Cal B(H_\alpha)$ and in this case all of our
invariants become trivial.  However, we point out that Powers
has recently worked out a {\it standard form} for \esg s
\cite{26} and any \esg\ in standard form admits an
invariant vector state.  Thus the case where $\Omega$ is
rank one should not be considered mysterious.
The other extreme case ($\Omega$ is of infinite rank)
appears to be amenable to the techniques developed below, but
there are technical difficulties associated with the infinite
rank case that require special attention.  We have elected to
postpone discussion of the infinite rank case to a subsequent
paper.

In broad terms, our results show that there is
an infinitesimal ``dynamical principle" that governs the behavior of
such a pair $(\alpha,\rho)$.  This invariant is analogous to the
Hamiltonian structure (on the cotangent bundle of a Riemannian
manifold) that governs the behavior of a constrained classical
mechanical system with a finite number of degrees of freedom.
Indeed, the work carried out in the following sections amounts to
a classification of Riemannian type structures in finite
dimensional matrix algebras.

The geometric nature of these invariants shows that the
``dynamical principles" that govern the behavior of \esg s
are closely akin to second order differential equations.  As
in classical mechanics, we are able to identify two
fundamental aspects of the dynamics: a metric term corresponding
to momentum and kinetic energy, and a potential term corresponding
to the driving force.
Moreover, one is free to specify the momentum space, the metric,
and the potential operator arbitrarily.  Regardless of how this is done
one obtains an \esg\ pair $(\alpha, \rho)$ (see
remark 8.8 for more detail).  In particular,
this gives an entirely new way of constructing \esg s.

On the other hand, not all pairs $(\alpha, \rho)$ (satisfying
the finiteness condition (0.4)) arise in this way.  But we
are able to give a useful characterization of those that do.
Broadly speaking, a pair $(\alpha, \rho)$ arises in this way iff
it satisfies two conditions.  First, $\alpha$ must obey a certain
{\it minimality} condition (which can always be arranged by
replacing $\alpha$ with a compression to a suitable hereditary
subalgebra of $\Cal B(H_\alpha)$).  More significantly,
$(\alpha,\rho)$ must also be {\it exact} in an appropriate sense.
Exactness is equivalent to several important properties, and it
is closely analogous to the hypothesis of classical mechanics in
which one supposes that the force (a vector field on configuration
space) should be the gradient of a potential.


\subheading{1.  Markov semigroups}
The purpose of this section is to make some observations
about the relationship that exists between \esg s and semigroups
of completely positive maps.

Let $\phi = \{\phi_t: t\geq 0\}$ be a semigroup of normal completely
positive linear maps of $\Cal B(H)$ such that $\phi_t(\bold 1) = \bold 1$,
and which is continuous in the sense that for every $a\in \Cal B(H)$,
$\xi, \eta\in H$, $ \<\phi_t(a)\xi,\eta\>$ is a continuous
function of $t$.  We will refer to such a semigroup simply as a
{\it completely positive semigroup}.

\proclaim{Definition 1.1}
A Markov semigroup is a pair $(\phi,\rho)$ consisting of a
completely positive semigroup $\phi$ acting on $\Cal B(H)$ and
a faithful normal state $\rho$ of $\Cal B(H)$ which is invariant
under $\phi$ in that $\rho\circ\phi_t = \rho$ for
every $t$.
\endproclaim

\remark{Remark}
We emphasize that the normal state should be {\it faithful}:
$\rho(a^*a) = 0 \implies a = 0$, for every $a\in\Cal B(H)$.
\endremark

Two Markov semigroups $(\phi,\rho)$, $(\phi^\prime,\rho^\prime)$
are said to be {\it conjugate} if there is a $*$-isomorphism
$\theta: \Cal B(H) \to \Cal B(H^\prime)$ such that
$$
\align
\theta\circ \phi_t &= \phi_t^\prime \circ \theta, \qquad t\geq 0,\tag{1.2.1}\\
\rho^\prime \circ \theta &= \rho.  \tag{1.2.2}
\endalign
$$
In this paper, we will be primarily concerned with
Markov semigroups acting on finite dimensional spaces; that is,
Markov semigroups which act on $n\times n$ matrix
algebras, $2\leq n < \infty$.

The connection between Markov
semigroups and \esg s is based on the following observations.
By an \esg\ we mean a semigroup $\alpha = \{\alpha_t: t\geq 0\}$
of normal $*$-endomorphisms of the algebra $\Cal B(H)$ of all
bounded operators on a separable Hilbert space $H$, such that
$\alpha_t(\bold 1) = \bold 1$ and which is continuous in the
sense described above.

\proclaim{Lemma 1.3} Let $\alpha$ be an \esg\ acting on
$\Cal B(H)$ and let $\omega$ be a normal state of $\Cal B(H)$
which is invariant under the action of
$\alpha_t$ for every $t\geq 0$.  Let
$p_0$ be the support projection of $\omega$.  Then
$\alpha_t(p_0) \geq p_0$ for every $t\geq 0$.
\endproclaim

\remark{Remark}
The support projection of $\omega$ is the smallest projection $p$
with the property that $\omega (\bold 1 - p) = 0$.  One has
$\omega(a) = \omega(p_0ap_0)$ for every $a\in \Cal B(H)$, and
$\omega(a^*a)=0$ iff $ap_0=0$.
\endremark

\demo{proof of 1.3}  Let $t\geq 0$, and consider the projection
$\alpha_t(p_0)$.  Because of the invariance of $\omega$ we have
$$
\omega(\bold 1 - \alpha_t(p_0)) = \omega\circ \alpha_t(\bold 1-p_0)
= \omega(\bold 1 - p_0) = 0,
$$
and hence $\alpha_t(p_0) \geq p_0$\qed
\enddemo

It follows immediately that the family of projections
$p_t = \alpha_t(p_0)$ is increasing on the interval
$0\leq t < \infty$.  Let us write $M$ for the von Neumann
algebra $\Cal B(H)$.  Then we can define a family
$\phi = \{\phi_t: t\geq 0\}$ of completely positive
maps
on the hereditary subalgebra $M_0 = p_0Mp_0$ as follows:
$$
\phi_t(a) = p_0\alpha_t(a)p_0, \qquad a\in M_0.
$$
We have $\phi_t(p_0)=p_0$, and the semigroup property follows
from the observation that if $s,t\geq 0$ then
$$
\phi_s(\phi_t(a)) = p_0\alpha_s(p_0\alpha_t(a)p_0)p_0
= p_0p_s\alpha_{s+t}(a)p_sp_0 = \phi_{s+t}(a),
$$
because $p_0p_s = p_sp_0 = p_0$.  Thus, if we let
$\rho$ be the restriction of $\omega$ to $M_0$ and then
identify $M_0$ with $\Cal B(p_0H)$, we find that we have
a Markov semigroup $(\phi,\rho)$.  If $\omega$ is weakly
continuous (equivalently, $p_0$ is finite dimensional) then
we may consider that $(\phi ,\rho)$ acts on an
$n\times n$ matrix algebra.

Notice that if $\omega$ happens to be a pure (vector) state,
then $p_0H$ is one-dimensional and $\phi$ is the trivial semigroup.
Thus all information has been lost in the passage from $(\alpha,\omega)$
to $(\phi,\rho)$.  On the other hand, in all other cases the Markov
semigroup will contain significant information about $(\alpha,\omega)$.

Indeed, it is often the case that $(\phi,\rho)$ completely determines
$(\alpha,\omega)$.  In order to discuss this
issue briefly, let $H_0 = p_0H$, and consider the subspace $H_+\subseteq H$
defined by
$$
H_+ = [\alpha_{t_1}(a_1)\alpha_{t_2}(a_2)\dots \alpha_{t_n}(a_n)\xi:
a_i\in M_0, \xi \in H_0, n = 1,2,\dots].
$$
There is a corresponding von Neumann subalgebra $M_+$ of $\Cal B(H)$,
namely the weakly closed algebra generated by the family of operators
$\{\alpha_t(a): a\in M_0, t\geq 0\}$.  It can be shown that $M_+$ is
a hereditary subalgebra of $\Cal B(H)$, and in fact
$$
M_+ = p_+\Cal B(H)p_+,
$$
$p_+$ denoting the projection of $H$ onto $H_+$.
The original pair $(\alpha, \omega)$ is called {\it minimal} if
$H_+ = H$, or equivalently $M_+ = \Cal B(H)$.

If $(\alpha,\omega)$ is minimal, then it is determined up to conjugacy
by its associated Markov semigroup $(\phi, \rho)$ (this follows from
work of Bhat described in the following paragraph).
But even when $(\alpha,\rho)$ is not
minimal, one may compress $\alpha$ to the invariant subalgebra
$M_+ = p_+\Cal B(H)p_+$ to obtain a minimal pair
$(\alpha^\prime,\omega^\prime)$.  Thus one may always assume that minimality
is satisfied provided one is willing to compress $\alpha$ to an
invariant corner of $\Cal B(H)$.

Moreover, a recent dilation theorem of B. V. R. Bhat \cite{6}, building
on and clarifying earlier partial results (\cite{17}, \cite{18}, \cite{19},
\cite{11}, \cite{14}, \cite{21}), implies that {\it every}
Markov semigroup $(\alpha,\rho)$ arises in this way from a minimal pair
$(\alpha,\omega)$, where $\alpha$ is an \esg\ and $\rho$ is a normal
$\alpha$-invariant state.  We will discuss Bhat's theorem, minimality
criteria, and properties of the ``$n$-point functions"
$$
\omega(\alpha_{t_1}(a_1)\alpha_{t_2}(a_2)\dots \alpha_{t_n}(a_n))
$$
defined for $a_1,\dots,a_n\in M_0$, $t_i\geq 0$, $n = 1,2,\dots$ in
a subsequent paper.

\subheading{2.  Differential operators on matrix algebras}

There is a natural notion of the {\it order} of a differential
operator which is easily adapted to the context of linear operators
on arbitrary commutative algebras.  Indeed, if $A$ is a complex commutative
algebra and $L:A\to A$ is a linear operator, one may define a sequence
of multilinear mappings $\Delta^nL:A^{n+1}\to A$ as follows:
$\Delta^0 L = L$ and
$$
\Delta^{n+1}L(f_1,\dots,f_n,f_{n+1};a) =
\Delta^nL(f_1,\dots,f_{n};f_{n+1}a) - f_{n+1}\Delta^nL(f_1,\dots,f_{n-1};a).
$$
One finds that $\Delta^nL(f_1,\dots,f_n;a)$ is a symmetric function
of its first $n$ variables $f_1,\dots,f_n$.  $L$ is said to be a differential
operator of order $n$ if $\Delta^{n+1}L = 0$ and $\Delta^{n}L \neq 0$.

The purpose of this section is to discuss the extent to which the
notion of order is appropriate for operators on noncommutative algebras,
and to introduce a suitable definition of {\it symbol} of a differential
operator.
For definiteness, we will consider the case where $A$ is the algebra
of all $n\times n$ matrices over $\Bbb C$, since it is this case that
is relevant for our purposes.  However, the reader will note that we
make no essential use of  structures specific to $M_n(\Bbb C)$.  Throughout,
$\Cal L(A)$ will denote the algebra of all linear operators
$L:A\to A$.

\proclaim{Definition} L is said to be a first order differential
operator if for every $a,x,y\in A$ we have
$$
L(xay)-xL(ay)-L(xa)y+xL(a)y = 0.  \tag{2.1}
$$
\endproclaim

Notice that the trilinear form appearing on the left side of
(2.1) is a noncommutative version
of the trilinear form $\Delta^2L(x,y;a)$ discussed above.  Moreover,
since $A$ has a unit (2.1) is equivalent to the somewhat
simpler condition
$$
L(xy) - xL(y) -L(x)y +xL(\bold 1)y = 0
$$
for all $x,y\in A$
(this will be discussed more fully later in the section).

A simple computation shows that any operator of the form
$$
L(x) = D(x) + ax \tag{2.2}
$$
where $a$ is a fixed element of $A$ and $D$ is a derivation must be
a first order differential operator.  The fact that the multiplier
$a$ appears on the left in (2.2) is inessential, since any operator
of the form $L(x) = D(x) + ax + xb$ can be put into the form (2.2) by
replacing $D$ with the derivation $D^\prime(x)  = D(x) + xb - bx$.
Noting that the first order differential
operators $L$ satisfying $L(\bold 1)=0$ are derivations,
it follows that
the first order differential operators are {\it exactly} those of the
form (2.2).

In contrast to the commutative case, operators having the form
$L(x) = ax$ or $L^\prime(x) = xa$ should  not be regarded
as ``order zero" differential operators,
since the difference of two such operators
$L(x) - L^\prime(x) = ax -xa$ is in this case a nontrivial derivation
which must be assigned order $1$.  Thus {\it there is no meaningful
definition of ``order zero" for operators on noncommutative algebras}.

Similarly, there is no viable concept of ``order $n$" for
$n > 2$.  For in the case where $A = M_n(\Bbb C)$ is a matrix
algebra, every linear operator $L$ on $A$
is a finite sum of double multipliers
$x\mapsto axb$, $a,b$ being fixed elements of $A$, and thus by
the preceding paragraph $L$ has the form
$$
L = A_1B_1 + A_2B_2+\dots + A_rB_r,
$$
where $A_k$ and $B_k$ are first order differential operators.  We
conclude that {\it every} operator on $A$ is of ``order" at most $2$.

Despite the somewhat negative tone of these remarks, we will find
Definition (2.1) to be quite useful.  In order to discuss this issue,
we recall the differential algebra $\Omega^*(A)$ associated
with $A$.  This was introduced in \cite{2}, and has
become a basic constituent of Connes' noncommutative differential
calculus \cite{9}.  Actually, we only require the two
modules $\Omega^1(A)$ and $\Omega^2(A)$, which can be defined for
unital $*$-algebras as follows.  Consider the tensor product
$A\otimes A$ as an involutive bimodule over $A$, with
$$
\align
a(x\otimes y)b &= ax\otimes yb,\\
(x\otimes y)^* &= y^*\otimes x^*.
\endalign
$$
The map $d:A\to A\otimes A$ defined by
$dx = \bold 1\otimes x - x\otimes \bold 1$ is a derivation for
which $(dx)^* = -d(x^*)$, and it is a universal derivation
of $A$ in the sense that if $E$ is any bimodule and $D:A\to E$ is
a derivation, then there is a unique homomorphism of bimodules
$\theta:\Omega^1(A)\to E$ such that $\theta\circ d = D$.  Every
element of $\Omega^1(A)$ is a finite sum of the form
$$
\omega = \sum_{k=1}^r a_k\,dx_k.
$$
Finally,
$\Omega^1(A)$ is the kernel of the multiplication map
$\mu:A\otimes A\to A$ defined by
$\mu(a\otimes b) = ab$, and thus we have an
exact sequence of involutive bimodules
$$
0 @>>> \Omega^1(A) @>>> A\otimes A @>>\mu> A @>>> 0.
\tag{2.3}
$$

For brevity, we will use the notation $\Omega^1$, $\Omega^2$ rather
than $\Omega^1(A)$, $\Omega^2(A)$ whenever it does not lead to confusion.
$\Omega^2$ is defined by
$$
\Omega^2 = \Omega^1\otimes_A \Omega^1,
$$
and every element of $\Omega^2$ is a finite sum of the form
$$
\omega = \sum_{k=1}^r a_k\,dx_k\,dy_k.
$$

Alternately, $\Omega^2$ can be viewed as the submodule of the
bimodule $A^{\otimes 3} = A\otimes A\otimes A$ that is linearly spanned
by elements of the form
$$
\align
adxdy &= a(\bold 1\otimes x-x\otimes \bold 1)(\bold 1\otimes y-y\otimes \bold
1)\\
&=a\otimes x\otimes y - ax\otimes \bold 1\otimes y
+ax\otimes y\otimes \bold 1 - a\otimes xy\otimes \bold 1.
\endalign
$$
There are two multiplication maps defined on $A^{\otimes 3}$ which
together give rise to a homomorphism of bimodules
$\mu: A^{\otimes 3} \to A^{\otimes 2}\oplus A^{\otimes 2}$:
$$
\mu(a\otimes b\otimes c) = (ab\otimes c,a\otimes bc),
$$
and we have an exact sequence of modules
$$
0 @>>> \Omega^2 @>>> A^{\otimes 3} @>>\mu>
A^{\otimes 2}\oplus A^{\otimes 2} @>>> 0.\tag2.4
$$
The involutions defined in $A^{\otimes 3}$ and
$A^{\otimes 2}\oplus A^{\otimes 2}$ by
$$
\align
(a\otimes b\otimes c)^* &= c^*\otimes b^*\otimes a^*\\
(\xi,\eta)^* &= (\eta^*,\xi^*)
\endalign
$$
make (2.4) into an exact sequence of involutive bimodules.

We will make frequent use of the following observation.
\proclaim{Lemma 2.5}
Let $L\in\Cal L(A)$.  There is a unique homomorphism of bimodules
$\theta_L\in \hom(\Omega^2,A)$ such that
$$
\theta_L(dx\,dy) = L(xy)-xL(y)-L(x)y+xL(\bold 1)y.
$$
\endproclaim

\demo{proof}
The uniqueness of $\theta_L$ is clear from the fact that $\Omega^2$
is spanned by elements of the form $a\,dx\,dy$, and
$$
\theta_L(a\,dx\,dy) = a(L(xy)-xL(y)-L(x)y+xL(\bold 1)y).
$$
For existence, we can define $\theta\in \hom(A^{\otimes 3},A)$ by
$\theta(a\otimes b\otimes c) = -aL(x)b$, and one finds that the restriction
of $\theta$ to elements of the form $dx\,dy$ is as required \qed
\enddemo

A straightforward computation shows that, more generally,
$$
\theta_L(a\,dx\,b\,dy\,c) =
aL(xby)c-axL(by)c-aL(xb)yc+axL(b)yc.
$$
Note too that  $L$ is a first order differential operator if, and only
if, $\theta_L=0$.

Now it follows from the definition of $\Omega^1$ that for
any derivation $D$
of $A$ there is a unique $\theta_D\in \hom(\Omega^1,A)$ such that
$$
\theta_D(dx) = D(x),\qquad x\in A.
$$
Moreover, the map $D\mapsto \theta_D$ is a linear isomorphism
of the space of derivations of $A$ onto $\hom(\Omega^1,A)$.  There is
a somewhat similar characterization of
$\hom(\Omega^2,A)$, valid for any unital
algebra $A$ for which the Hochschild cohomology space $H^2(A,A)$ is trivial:

\proclaim{Lemma 2.6}
Let $\theta\in \hom(\Omega^2,A)$.  There is a linear operator
$L\in \Cal L(A)$ satisfying $L(\bold 1)=0$ and $\theta = \theta_L$.
L is unique up to a perturbation of the form $L^\prime = L+D$ where
$D$ is a derivation of $A$.
\endproclaim
\demo{proof}
Fix $\theta$, and consider the bilinear map $T:A\times A\to A$ defined by
$T(x,y) = \theta(dx\,dy)$.  The Hochschild coboundary of $T$
$$
\align
bT(x,y,z) &= xT(y,z) -T(xy,z)+T(x,yz)-T(x,y)z\\
&=\theta(x\,dy\,dz-d(xy)\,dz +dx\,d(yz)-dx\,dy\,z)
\endalign
$$
vanishes because
$$
\align
d(xy)\,dz &= x\,dy\,dz+dx\,y\,dz,\qquad \text{and}\\
dx\,d(yz) &= dx\,y\,dz + dx\,dy\,z.
\endalign
$$
Since $H^2(A,A)=0$ when $A$ is a matrix algebra, there exists
a linear operator $L\in \Cal L(A)$ for which
$$
\theta(dx\,dy) = bL(x,y) = xL(y)-L(xy)+L(x)y.
$$
Setting $L_0(x) = L(x)-xL(\bold 1)$ we find that $L_0(\bold 1) = 0$
and $\theta = \theta_{-L_0}$.

If $L_1$ and $L_2$ are two linear operators satisfying $L_1(\bold 1) =
L_2(\bold 1) = 0$ and $\theta = \theta_{L_1}=\theta_{L_2}$, then
$\theta_{L_1-L_2}=0$.  Since $L_1-L_2$ is a first order differential
operator that annihilates $\bold 1$, it must be a derivation\qed
\enddemo

These remarks about the differential nature of linear operators on
$A$ are summarized in the following exact sequence of complex vector
spaces, in which $\Cal L_0(A) = \{L\in \Cal L(A): L(\bold 1) =0\}$,
$\iota$ is the identification of $\hom(\Omega^1,A)$ with
derivations by way of $\iota\theta(x)=\theta(dx)$ $x\in A$,
and $\theta:L\mapsto \theta_L$ is the mapping of Lemma 2.5,
$$
0 @>>> \hom(\Omega^1,A)@>>\iota> \Cal L_0(A)
@>>\theta> \hom(\Omega^2,A)@>{}>> 0.
$$

Throughout the remainder of this paper we will be concerned with
{\it pairs} $(A,\rho)$ consisting of a finite dimensional \cstar\
$A$ (usually $M_n(\Bbb C)$) together with
a distinguished {\it faithful} state $\rho$, that is a linear
functional satisfying
$$
\align
\rho(\bold 1) &= 1,\qquad \text{and}\\
\rho(a^*a)&>0,\qquad \text{for every nonzero }a\in A.
\endalign
$$
The associated space $\Cal D(A,\rho)$ of differential operators
is defined by the three requirements:
$$
\alignat 2
&\text{Normalization:} &\qquad L(\bold 1) &= 0 \tag{2.7.1}\\
&\text{Divergence zero:} &\qquad\rho\circ L&=0 \tag{2.7.2}\\
&\text{Symmetry:} &\qquad L(x^*) &= L(x)^*,\quad x\in A.   \tag{2.7.3}
\endalignat
$$
Conditions (2.7.1) and (2.7.2) are dual to each other, in a sense
that will be exploited in the following sections.  In general,
an operator $L$ satisfying (2.7.3) will be called {\it symmetric}.
$\Cal D(A,\rho)$ is a real vector space of linear operators on $A$.
For every $L\in \Cal D(A,\rho)$ we define the {\it symbol} of
$L$ by $\sigma_L = \rho\circ \theta_L$.  More explicitly, $\sigma_L$
is the linear functional defined on $\Omega^2$ by
$$
\sigma_L(a\,dx\,dy) = \rho(aL(xy)-axL(y)-aL(x)y).
$$
We collect the following elementary properties of the symbol for
later use.

\proclaim{Proposition 2.9}
For every $L\in \Cal D(A,\rho)$, $\sigma_L$ has the following properties:
$$
\align
\sigma_L(\xi^*) &= \overline{\sigma_L(\xi)}, \qquad \xi\in \Omega^2\tag{i}\\
\sigma_L&=0 \text{ iff $L$  is a derivation}.\tag{ii}
\endalign
$$
\endproclaim
\demo{proof}
(i) follows from the formula $\theta_L(\xi^*) = \theta_L(\xi)^*$,
which in turn reduces to the property $L(x^*) = L(x)^*$ after evaluating
both sides at elements of the form $\xi = dx\,dy$.

For (ii), we note that the bilinear form $a,b\in A\mapsto \rho(ab)$
is nondegenerate because $\rho$ is a faithful state.  Thus
$\theta_L(dx\,dy) = 0$ for all $x,y\in A$ iff
$\rho(a\theta_L(dx\,dy))= 0 $ for every $a,x,y\in A$, hence (ii).\qed
\enddemo

Finally, we note that for any finite set $D_0,D_1,\dots,D_r$ of
symmetric derivations of $A$ satisfying $\rho\circ D_k = 0$ for
$0\leq k\leq r$, the operator
$$
L(x) = D_0(x) + \sum_{k=1}^rD_k^2(x)
$$
belongs to $\Cal D(A,\rho)$, and has symbol
$$
\sigma_L(a\,dx\,dy) = 2\sum_{k=1}^r \rho(aD_k(x)D_k(y)).  \tag{2.10}
$$
More generally, one has
$$
\sigma_L(a\,dx\,b\,dy\,c)= 2\sum_{k=1}^r \rho(aD_k(x)bD_k(y)c).  \tag{2.11}
$$
It is operators of this type that will be central to our analysis
of Markov semigroups.

\subheading{3.  Elliptic operators}
Throughout this section, $A$ will denote a
finite dimensional matrix algebra and $(\phi,\rho)$ will
denote a Markov semigroup
acting on $A$.  In this case
$\phi$ is uniformly continuous
$$
\lim_{t\to 0}\|\phi_t - \text{id}\| = 0.
$$
The infinitesimal generator $L$ exists in all senses, and obeys
the three properties (2.7).  Hence $L$ belongs to $\Cal D(A,\rho)$.
There are two characterizations of the generators of completely
positive semigroups that are significant for our purposes.  The first is due
to Lindblad \cite{20} and independently to Gorini {\it et al}
\cite{15} (also see \cite{10, Theorem 4.2}).
The second characterization is due to Evans and Lewis \cite{14},
based on work of Evans \cite{12}. These two
results can be paraphrased in our context as follows.

\proclaim{Theorem 3.1}An operator $L\in \Cal D(A,\rho)$
generates a Markov semigroup iff there is a completely
positive linear map $P: A\to A$ and an operator $a\in A$ such that
$$
L(x) = P(x) + ax + xa^*,\qquad x\in A.
$$
\endproclaim

\proclaim{Theorem 3.2}An operator $L\in \Cal D(A,\rho)$ generates a
Markov semigroup iff for every $n\geq1$ and every set of
elements $a_1,b_1,\dots, a_n,b_n\in A$,
satisfying $b_1a_1+b_2a_2+\dots+b_na_n=0$ we have
$$
\sum_{i,j=1}^n a_j^*L(b_j^*b_i)a_i \geq 0.
$$
\endproclaim

An operator $L$ satisfying the condition of (3.2) is called {\it conditionally
completely positive} \cite{13}.  While 3.1 tells us exactly which operators
are generators of Markov semigroups, the cited decomposition of $L$ into
a sum of more familiar operators is unfortunately not unique.

The purpose of
this section is to use (3.2) to give a new characterization of generators
of Markov semigroups in terms of their symbols $\sigma_L:\Omega^2(A)\to \Bbb
C$.
Recall that the involution in $\Omega^1$ is defined by
$$
(a\,dx)^* = -dx^*\,a^* = -d(x^*a^*) + x^*\,da^*,
$$
while that of $\Omega^2$ is defined by
$$
(a\,dx\,dy)^* = dy^*\,dx^*\,a^*.
$$
Thus for any two 1-forms $\omega_2, \omega_2\in \Omega^1$ we can
form various products in $\Omega^2$: $\omega_1\omega_2$, $\omega_1^*\omega_2$,
etc.

\proclaim{Theorem 3.3}  An operator $L\in D(A,\rho)$ generates
a Markov semigroup on $(A,\rho)$ iff $\sigma_L(\omega^*\omega)\leq 0$
for every $\omega\in \Omega^1$.
\endproclaim

\demo{proof}Let $L$ be an arbitrary operator in $\Cal L(A)$,
and let $\theta_L:\Omega^2\to A$ be the module homomorphism defined
in Lemma 2.5.  Choose any sequence of elements
$a_1,b_1,\dots,a_n,b_n\in A$ satisfying
$\sum_kb_ka_k=0$ and define an element $\omega\in A\otimes A$
by
$$
\omega = \sum_{k=1}^n b_k\otimes a_k.  \tag{3.5}
$$
Notice that $\omega$ belongs to the kernel of the multiplication map
$\mu:A\otimes A\to A$ and hence $\omega\in \Omega^1$.  Conversely,
every element $\omega\in \Omega^1$ can be deomposed into a sum
of the form
(3.5) which belongs to the kernel of the multiplication map.
Now for such a 1-form $\omega$ we have
$$
\omega^*\omega = \sum_{i,j=1}^na_j^*\otimes b_j^*b_i\otimes a_i.
$$
Let $\theta\in \hom(\Omega^2, A)$ be defined by
$\theta(a\otimes x\otimes b) = -aL(x)b$.  The proof of Lemma
2.5 shows that $\theta_L$ is obtained by restricting $\theta$
to $\Omega^2$, hence
$$
\theta_L(\omega^*\omega) = -\sum_{i,j=1}^n a_j^*L(b_j^*b_i)a_i.
$$

This observation, together with Theorem 3.2, shows that an operator
$L\in \Cal D(A,\rho)$ generates a Markov semigroup if, and only if,
$\theta_L(\omega^*\omega) \leq 0$ for every $\omega\in \Omega^1$.
It follows that for every generator $L$ and every $\omega\in \Omega^1$
we have $\sigma_L(\omega^*\omega) = \rho(\theta_L(\omega^*\omega)\leq 0$.

Conversely, if $\sigma_L(\omega^*\omega)\leq 0$ for every
$\omega\in \Omega^1$ then for every $b\in A$ we have
$$
\rho(b^*\theta_L(\omega^*\omega)b) =
\rho(\theta_L((\omega b)^*\omega b)) = \sigma_L((\omega b)^*\omega b)\leq 0.
$$
Since $\rho$ is a faithful state and the operator $T=\theta_L(\omega^*\omega)$
obeys $\rho(b^*Tb) \leq 0$ for every $b$, it follows that $T\leq 0$.
Hence $\theta_L(\omega^*\omega) \leq 0$ for every $\omega\in \Omega^1$ and
we may conclude from the preceding paragraphs that $L$ generates a Markov
semigroup \qed
\enddemo

In view of Theorem 3.3, we make the following
\proclaim{Definition 3.4} An {\rm elliptic} operator
is an operator $L\in \Cal D(A,\rho)$ satisfying
$$
\sigma_L(\omega^*\omega)\leq 0,\qquad \text{for every }\omega\in \Omega^1.
$$
\endproclaim

The classification problem for Markov semigroups is now reduced to
the problem of classifying elliptic operators up to the natural notion
of conjugacy: $L\in \Cal D(A,\rho)$ and
$L^\prime\in \Cal D(A^\prime,\rho^\prime)$ are said to be {\it conjugate}
if there is a $*$-isomorphism $\theta:A\to A^\prime$ satisfying
$$
\align
\rho^\prime\circ \theta &= \rho\\
L^\prime\circ \theta &= \theta\circ L.  \tag{3.5}
\endalign
$$

\subheading{4.  Momentum}
Let $(A,\rho)$ be a matrix algebra endowed with a faithful
state $\rho$.  In this section we show how elliptic operators
are constructed from more basic structures.

Let $p$ be an element of $A$ satisfying $p^*=-p$.  If we
write $D_p(a) = [p,a] = pa-ap$, then $D_p$ is a symmetric
derivation of $A$; replacing $p$ with $p-\rho(p)\bold 1$ if
necessary, we can assume that $p$ is normalized so that
$\rho(p)=0$.  With this convention for normalization,
the operator $p$ is uniquely determined by
the derivation $D_p$.  Letting $\tau$ be the tracial
state of $A$, we can define the density matrix $h$ of
$\rho$ by $\rho(a) = \tau(ha)$, $a\in A$.  $h$ is a
self adjoint matrix with strictly positive spectrum.
Notice, finally, that $\rho\circ D_p = 0$ iff $p$
commutes with $h$; and in that case we have
$D_p\in \Cal D(A,\rho)$.

\proclaim{Definition 4.1}  A momentum space is a pair
$(P,\<,\>)$ consisting of a real linear space $P$ of
skew-adjoint operators $p\in A$ satisfying $\rho(p)=0$
and $\rho\circ D_p = 0$, together with a real inner product
$\<,\>:P\times P\to \Bbb R$.
\endproclaim

We emphasize that the inner product on $P$ can be specified
arbitrarily, and in particular if $A$ is realized concretely
as $\Cal B(H)$ for a finite dimensional Hilbert space $H$ then
there need be no relation between the inner products on $P$ and
$H$.  Given a momentum space $(P,\<,\>)$, we can construct an
operator $\Delta\in \Cal D(A,\rho)$ as follows.  Choose an
orthonormal basis $p_1,p_2,\dots,p_r$ for $P$ and let
$D_k(a) = [p_k,a]$, $1\leq k\leq r$.  $\Delta$ is defined by
$$
\Delta = \sum_{k=1}^r D_k^2.  \tag{4.2}
$$

\proclaim{Proposition 4.3}
The operator $\Delta$ of {\rm(4.2)} does not depend on the
orthonormal basis chosen for $P$, and its symbol obeys
$$
\sigma_\Delta(dx\,dy) = -2\rho(\Delta(x) y),\tag{4.3.1}
$$
for all $x,y\in A$.
\endproclaim

\demo{proof}
Since $D_k$ is a derivation in $\Cal D(A,\rho)$ we have
$$
\rho(D_k^2(x)y) = \rho(D_k(D_k(x)y)) - \rho(D_k(x)D_k(y))
= -\rho(D_k(x)D_k(y))
$$
for each $k$, so that by (2.10)
$$
\rho(\Delta(x)y) = -\sum_{k-1}^r \rho(D_k(x)D_k(y)) =
-\frac{1}{2}\sigma_\Delta(dx\,dy).
$$
This establishes (4.3.1).

Let $p^\prime_1,\dots,p^\prime_r$ be another orthonormal basis
for $P$ and set $\Delta^\prime = \sum_kD_k^{\prime 2}$, where
$D_k^\prime = [p_k^\prime,\cdot]$.  Since $\rho$ is a faithful state
it suffices to show that $\rho(\Delta(x)y) = \rho(\Delta^\prime(x)y)$
for all $x,y\in A$; and by (4.3.1) this will follow from the assertion
$$
\sum_{k=1}^r\rho(D_k^\prime(x)D_k^\prime(y)) =
\sum_{k=1}^r\rho(D_k(x)D_k(y)),
$$
for all $x,y\in A$.
To prove the latter we make the substitution
$p_k^\prime = \sum_i\<p_k^\prime,p_i\> p_i$ in the expression
$D^\prime_k(x)D_k^\prime(y) = [p_k^\prime,x][p_k^\prime,y]$ to obtain
$$
D_k^\prime(x)D_k^\prime(y) = \sum_{i,j=1}^r\<p_k^\prime,p_i\>
\<p_k^\prime,p_j\>D_i(x)D_j(y).
$$
When the right side is summed on $k$, we may use
the orthonormality of $\{p_k\}$ and $\{p_k^\prime\}$ in the form
$\sum_k\<p_k^\prime,p_i\> \<p_k^\prime,p_j\> = \delta_{ij}$ to
obtain
$$
\sum_{k=1}^rD_k^\prime(x)D_k^\prime(y)= \sum_{k=1}^rD_k(x)D_k(y),
$$
and the claim follows\qed
\enddemo

The operator $\Delta$ of equation (4.2) is called the {\it Laplacian}
of the momentum space $(P,\<,\>)$.  Finally, let $v$ be any
skew-adjoint element of $A$ for which $D_v(a) = [v,a]$ obeys
$\rho\circ D_v = 0$, and set
$$
L(x) = \Delta(x) + [v,x],\qquad x\in A.\tag{4.3}
$$

\proclaim{Proposition 4.4}  The operator $L$ of \rm{(4.3)} is an
elliptic operator in $\Cal D(A,\rho)$.
\endproclaim
\demo{proof}
It is obvious that $L$ satisfies the criteria (2.7) for membership
in $\Cal D(A,\rho)$.

Since $D_v$ is a derivation we have $\sigma_L = \sigma_\Delta$, and
thus it suffices to show that $\Delta$ is elliptic.  For that,
choose
$$
\omega = \sum_{k=1}^n a_k\,dx_k
$$
in $\Omega^1$.  Then
$$
\omega^*\omega = \sum_{i,j=1}^n(dx_i)^*\,a_i^*a_j\,dx_j
= -\sum_{i,j=1}^n dx_i^*\,a_i^*a_j\,dx_j.
$$
Thus by (2.11) we have
$$
\sigma_L(\omega^*\omega) = -\sum_{k=1}^r\sum_{i,j=1}^n
\rho(D_k(x_i^*)a_i^*a_jD_k(x_j)) = -\sum_{k=1}^r \rho(z_k^*z_k)\leq 0,
$$
where $z_k = \sum_ja_jD_k(x_j)$, using the fact
that $D_k(x^*) = D_k(x)^*$\qed
\enddemo

The following observation shows that an operator of the form
(4.3) determines both of its summands uniquely.  $L^2(A,\rho)$ denotes
the finite dimensional Hilbert space obtained by endowing  $A$ with
the inner product $\<a,b\>_\rho = \rho(b^*a)$.

\proclaim{Proposition 4.5}
Let $L = \Delta + [v,\cdot]$ have the form \rm{(4.3)}, and consider
$L$ as an operator on the Hilbert space $L^2(A,\rho)$.  Then
$L+L^* = 2\Delta$ and $L-L^* = 2[v,\cdot]$.
\endproclaim
\demo{proof}
It suffices to show that
$\<\Delta(a),b\>_\rho = \<a,\Delta(b)\>_\rho$ and
$\<[v,a],b\>_\rho = -\<a,[v,b]\>_\rho$.  But if $z$ is any
element of $A$ satisfying $z^*=-z$ and $\rho([z,a]) = 0$ for
all $a\in A$, then the operator $D(a) = [z,a]$ is a symmetric
derivation which induces a skew-adjoint operator in $\Cal B(L^2(A,\rho))$.
Indeed,
$$
\<D(a),b\>_\rho = \rho(b^*D(a)) = \rho(D(b^*a))-\rho(D(b^*)a)
=0-\rho(D(b)^*a)=-\<a,D(b)\>_\rho.
$$
Therefore $D^2$ is a self-adjoint operator on $L^2(A,\rho)$.
It follows that $\Delta$ is a self-adjoint operator on $L^2(A,\rho)$
and $[v,\cdot]$ is skew-adjoint.  (4.5) follows\qed
\enddemo

\subheading{5. Modular Cohomology}
The purpose of this section is to discuss certain
cohomological issues so as to provide a context
for the following section.

\proclaim{Definition 5.1}
A modular algebra is a pair $(A,\rho)$ consisting of a unital
$*$-algebra $A$ and a faithful state $\rho$ on $A$ with the
following property:
for every element $a\in A$ there is an element $\delta(a)\in A$
such that
$$
\rho(ab) = \rho(b\delta(a)),\qquad b\in A.  \tag{5.1.1}
$$
\endproclaim

\remark{Remarks}
By a faithful state on $A$ we mean a linear functional $\rho$
satisfying $\rho(\bold 1) = 1$ and $\rho(a^*a)>0$ for every
$a\neq 0$ $\in A$.  Since $\rho$ is faithful the bilinear form
$$
a,b\mapsto \rho(ab)
$$
is nondegenerate, and hence the element $\delta(a)$ defined by
(5.1.1) is {\it unique}.
\endremark

Straightforward calculations show that $\delta$ is an automorphism
of the algebra structure of $A$ for which $\rho\circ\delta = \rho$,
and from the property $\rho(x^*) = \overline{\rho(x)}$ one readily
deduces
$$
\delta(a)^* = \delta^{-1}(a^*), \tag{5.2}
$$
$\delta^{-1}$ denoting the inverse automorphism.  $\delta$ is
called the {\it modular automorphism} of $(A,\rho)$.  We deviate
from the traditional notation $\Delta$ for the modular automorphism
associated with a faithful normal state of a von Neumann algebra
in order to reserve $\Delta$ for the Laplacian of a momentum space.

\example{Example 5.3}
For our immediate purposes we have the case where $A$ is the
$*$-algebra $M_n(\Bbb C)$ of $n\times n$ complex matrices and
$\rho$ is a faithful state.  If $\tau$ is the normalized trace
on $M_n(\Bbb C)$ then we can define an operator $h\in M_n(\Bbb C)$
by $\rho(a) = \tau(ha)$, $a\in A$.  $h$ is a self-adjoint matrix with
strictly positive spectrum and $\delta$ is the inner automorphism
$\delta(a) = h a h^{-1}$.
\endexample

\example{Example 5.4}
Let $M$ be an arbitrary von Neumann algebra
and let $\rho$ be a faithful normal
state of $M$.  Let $\sigma^\rho$ be the modular automorphism group
of $\rho$.  We define $A$ to be the set of all elements $a\in M$
whose spectrum relative to this group is compact; equivalently, $A$
consists of all elements $a$ for which there is a function
$f = f_a\in L^1(\Bbb R)$ whose Fourier transform has compact
support and for which
$$
\int_{-\infty}^{+\infty} f(t)\,\sigma^\rho_t(a)\,dt = 0.
$$
$A$ is a unital $*$-algebra for which $\sigma^\rho_t(A) = A$ for
every $t\in \Bbb R$.  For each $a\in A$, the function
$t\in \Bbb R\mapsto \sigma^\rho_t(a)$ extends uniquely to an
entire function
$$
z\in \Bbb C\mapsto \sigma_z^\rho(a)\in A.
$$
$\{\sigma_z^\rho: z\in \Bbb C\}$ defines a group of automorphisms of
$A$ parameterized by the additive group of complex numbers, satisfying
$\sigma_z^\rho(a)^* = \sigma_{\overline z}^\rho(a^*)$.  The KMS condition
as formulated in \cite{23, 8.12.2} implies that the automorphism
$\delta = \sigma_i$ satisfies
$$
\rho(ab) = \rho(b\delta(a)),\qquad a,b\in A.
$$
Therefore $(A,\rho)$ is a modular algebra.
\endexample

Example 5.4 shows that the considerations of this section
are appropriate in situations where there is no trace whatsoever,
such as that in which one has an \esg\ acting on a factor of
type $III$.  We intend to take up these more general issues in
a subsequent paper.

We first sketch a generalization of Connes' cyclic cohomology that
is appropriate for modular algebras $(A,\rho)$.  By a
{\it modular cochain} of dimension $n\geq 0$ we mean a multilinear
functional of $n+1$ variables $\phi:A^{n+1}\to \Bbb C$ satisfying
$$
\phi(a^0,a^1,\dots,a^n) = (-1)^n\phi(\delta^{-1}(a^n),a^0,\dots,a^{n-1}).
\tag{5.5}
$$
The formula (5.5) is ambiguous in the case $n=0$, and in that case we
intend that $\phi$ should satisfy $\phi(a) = \phi(\delta^{-1}(a))$.
$C^n(A)$ denotes the vector space of $n$-dimensional cochains.  For
every $\phi\in C^n(A)$ we define a coboundary $b_\delta\phi$ by
$$
\align
b_\delta\phi(a^0,\dots,a^{n+1}) = &
\sum_{k=0}^n (-1)^k\phi(a^0,\dots,a^ka^{k+1},\dots,a^{n+1})  \tag{5.6} \\
&+ (-1)^{n+1}\phi(\delta^{-1}(a^{n+1})a^0,a^1,\dots,a^n).
\endalign
$$
One can show that the operator $b_\delta$ maps $C^n(A)$ into
$C^{n+1}(A)$ and satisfies $b_\delta^2 = 0$.  The cohomology
of the resulting complex is called {\it modular cohomology}.
A multilinear form $\phi: A^{n+1}\to \Bbb C$ satisfying (5.5)
and $b_\delta\phi = 0$ is
called a modular cocycle.  $\phi$ is called {\it exact} if there
is an $n-1$ cochain $\psi\in C^{n-1}(A)$ such that $\phi = b_\delta\psi$.

Notice that if $\rho$ is a trace on $A$ then $\delta$ is the identity
automorphism and we have cyclic cohomology (\cite{9, pp. 182--190}).
In the case where $\rho$ is a faithful state on a matrix algebra the
resulting modular cohomology is trivial; nevertheless these considerations
(especially for dimensions zero, one and two) will be relevant for our
purposes below.

A zero cochain is a linear functional $\phi\in A^\prime$ satisfying
$\phi\circ\delta = \phi$.  Because $\rho$ is faithful there is a unique
operator $p\in A$ for which $\phi(a) = \rho(ap)$, and the condition
$\phi\circ \delta = \phi$ is equivalent to $\delta(p) = p$; i.e.,
$p$ commutes with the ``density matrix" of $\rho$.

A zero cocycle is simply a scalar multiple of $\rho$,
$\phi(a) = \lambda\rho(a)$ for some $\lambda \in \Bbb C$.
We will require the following simple characterization of
exactness for one-cochains.

\proclaim{Proposition 5.7}
Let $\phi \in C^1(A)$.  Then $\phi$ is exact iff it has the form
$$
\phi(a,x) = \rho(a[v,x]),\tag{5.7.1}
$$
where $v$ is an element of $A$ satisfying $\rho([v,a]) = 0$ for
every $a\in A$.
\endproclaim

\demo{proof}
Supose first that $p\in A$ satisfies $\rho([p,a])=0$ for all
$a\in A$ (equivalently, $p$ commutes with the density matrix of
$\rho$) and consider the form $\phi(a,x) = \rho(a[p,x])$.
Let $\psi$ be the linear functional $\psi(z) = -\rho(zp)$.
Since $\delta(p)=p$ we have
$$
\psi(\delta(z))=-\rho(\delta(z)p)
=-\rho\circ\delta(zp)=-\rho(zp)=\psi(z),
$$
and hence $\psi\in C^0(A)$.  The coboundary of $\psi$ is given by
$$
\align
b_\delta\psi(a,x) &= \psi(ax)-\psi(\delta^{-1}(x)a)=
-\rho(axp)+\rho(\delta^{-1}(x)ap)\\
&=\rho(apx)-\rho(axp) = \rho(a[p,x]),
\endalign
$$
and hence $\phi = b_\delta\psi$.

Conversely, suppose $\phi$ is exact and choose $\psi\in C^0(A)$ such that
$\phi(a,x) = \psi(ax) - \psi(\delta^{-1}(x)a)$.  Since $\rho$ is faithful
there is an element $p\in A$ for which $\psi(z) = -\rho(zp)$.  One has
$\delta(p) = p$, and the preceding argument can be reversed to obtain
formula (5.7.1)\qed
\enddemo

\subheading{6. Exactness}
A central issue in this work has been to appropriately
characterize the elliptic operators $L\in \Cal D(A,\rho)$
that can be decomposed as in (4.3)
$$
L = \Delta + [v,\cdot]
$$
as a first order perturbation of the Laplacian
of some momentum space $(P,\<\cdot,\cdot\>)$.  In this
section we introduce the notion of {\it exact} differential
operator and give several characterizations of exactness
(Theorem 6.9).  In section 7 we show that the exact elliptic
operators are precisely those that can be so decomposed.

Throughout this section, $(A,\rho)$ will denote a matrix
algebra $A$ endowed with a faithful state $\rho$.  Let
$L\in \Cal D(A,\rho)$ be a differential operator.  We define
$\omega_L$ as the following bilinear form on $A$,
$$
\omega_L(x,y) = \rho(xL(y)) - \rho(L(x)y),\qquad x,y\in A.
$$
$\omega_L$ will occupy a central role in the remainder of
this paper.

Any operator $L\in \Cal L(A)$ can be regarded as an operator
on the Hilbert space $L^2(A,\rho)$ and as such it has an
adjoint $L^*\in \Cal B(L^2(A,\rho))$.  If $L$ is a {\it symmetric}
operator (i.e., $L(x^*) = L(x)^*$ for every $x\in A$) then
$L^*$ is uniquely defined by
$$
\rho(L(x)y) = \rho(xL^*(y)),\qquad x,y\in A.  \tag{6.1}
$$
Notice that (6.1) is false if $L$ is not symmetric, but
it can be repaired in that case by replacing $L(x)$ on the
left side with $L(x^*)^*$.

We are interested in operators that belong to $\Cal D(A,\rho)$.
For such an $L$ we find that $L^*(\bold 1) = 0$ by setting
$y=\bold 1$ in (6.1) and using $\rho\circ L = 0$.  Similarly,
$\rho\circ L^* = 0$ follows from the condition $L(\bold 1) = 0$.
However, the adjoint of an operator in $\Cal D(A,\rho)$ need
{\it not} symmetric.  That information is contained in
the form $\omega_L$ as follows.

\proclaim{Proposition 6.2}
Let $L\in \Cal D(A,\rho)$.  Then $L^*\in \Cal D(A,\rho)$ if and
only if $\omega_L$ is a modular one-cochain.
\endproclaim
\demo{proof}
The condition on $\omega_L$ means
$$
\omega_L(x,y) = -\omega_L(\delta^{-1}(y),x),\qquad x,y\in A.\tag{6.3}
$$
The left side of (6.3) is
$\rho(xL(y)) - \rho(L(x)y)$ while the right side is
$$
\align
-\rho(\delta^{-1}(y)L(x))+\rho(L(\delta^{-1}(y))x) &=
-\rho(L(x)y)+\rho(\delta^{-1}(y)L^*(x))\\
&=-\rho(L(x)y)+\rho(L^*(x)y).
\endalign
$$
Setting these two expressions equal we find that (6.3) is
equivalent to
$$
\rho(xL(y)) = \rho(L^*(x)y).\tag{6.4}
$$
The left side of (6.4) can be rewritten as follows,
$$
\align
\rho(xL(y)) &= \overline{\rho(L(y)^*x^*)} = \overline{\<x^*,L(y)\>}_\rho
=\overline{\<L^*(x^*),y\>}_\rho \\
&= \overline{\rho(y^*L^*(x^*))} = \rho(L^*(x^*)^*y).
\endalign
$$
Thus (6.4) is equivalent to the formula $L^*(x^*)^* = L^*(x)$,
i.e., that $L^*$ should be a symmetric operator \qed
\enddemo

\proclaim{Definition 6.5}
An operator $L\in \Cal D(A,\rho)$ is called {\rm exact} if there
is a modular zero-cochain $\phi\in C^0(A)$ such that
$\omega_L = b_\delta \phi$.
\endproclaim

\remark{Remarks}
Since the coboundary operator $b_\delta$ maps cochains to
cochains, we see from Proposition 6.2 that for any exact operator
$L\in \Cal D(A,\rho)$ we have $L^*\in \Cal D(A,\rho)$.
Exactness is considerably stronger, however, and the remainder
of this section is devoted to giving several more concrete
characterizations of exact differential operators.
\endremark

We have already pointed out that the modular cohomology of a
matrix algebra is trivial.  Hence $L$ {\it is exact iff }
$\omega_L$ {\it is a modular one-cocycle}.   Equivalently,
an operator $L\in \Cal D(A,\rho)$ is exact iff the following
two conditions are satisfied
$$
\align
&L^*\in \Cal D(A,\rho) \tag{6.6.1}\\
&b_\delta\omega_L = 0.\tag{6.6.2}
\endalign
$$

Let $\delta$ be the modular automorphism of $(A,\rho)$.
$\delta$ induces a natural linear isomorphism
$\hat\delta:\Omega^1\to \Omega^1$ defined by
$$
\hat\delta(\sum_ja_j\,dx_j) = \sum_j\delta(a_j)\,d\delta(x_j).\tag{6.7}
$$
One has $\hat\delta(a\omega b) = \delta(a)\hat\delta(\omega)\delta(b)$
for $a,b\in A$, $\omega\in \Omega^1$.  Notice too that the
formula $\delta(a)^* = \delta^{-1}(a^*)$ for $a\in A$ implies that
$$
\hat\delta(\omega^*) = \hat\delta^{-1}(\omega)^*,\qquad \omega\in
\Omega^1,\tag{6.8}
$$
where $\hat\delta^{-1}$ in (6.8) denotes the inverse of
$\hat\delta:\Omega^1\to\Omega^1$, defined by
$\hat\delta^{-1}(a\,dx) = \delta^{-1}(a)\,d\delta^{-1}(x)$,
$a,x\in A$.

\proclaim{Theorem 6.9}
For any differential operator $L$ in $\Cal D(A,\rho)$, the following
are equivalent.

\roster
\item"{(i)}"
$L$ is exact.

\item"{(ii)}"
There is an element $v\in A$ satisfying $v^*=-v$, $\delta(v) = v$,
and
$$
\omega_L(x,y) = \rho(x[v,y]),\qquad x,y\in A.
$$
\item"{(iii)}"
$L-L^*$ is a derivation.
\item"{(iv)}"The symbol of $L$ satisfies the KMS condition
$$
\sigma_L(\omega_1\omega_2) = \sigma_L(\omega_2\hat\delta(\omega_1)),
\qquad \omega_1,\omega_2\in \Omega^1.
$$
\endroster
\endproclaim

For the proof, we require the following two results about symbols.

\proclaim{Lemma 6.10}
Let $L\in \Cal L(A)$ be an arbitrary operator.  The symbol
of $L$ obeys
$$
\sigma_L(a\xi) = \sigma_L(\xi\,\delta(a)),\qquad a\in A,\xi\in \Omega^2.
$$
\endproclaim
\demo{proof}
Using the homomorphism of modules $\theta_L\in \hom(\Omega^2,A)$ of
Lemma 2.5 and the fact that $\sigma_L = \rho\circ\theta_L$ we have
$$
\align
\sigma_L(a\xi) &= \rho(\theta_L(a\xi)) =\rho(a\theta_L(\xi))
=\rho(\theta_L(\xi)\delta(a)) \\
&= \rho(\theta_L(\xi \,\delta(a))=\sigma_L(\xi\,\delta(a)),
\endalign
$$
as asserted\qed
\enddemo

\proclaim{Lemma 6.11}
Let $L$ be an operator in $\Cal D(A,\rho)$ and let
$L^*\in \Cal B(L^2(A,\rho))$ be its adjoint.  Then we have
$$
\sigma_L(\omega_1\omega_2) = \sigma_{L^*}(\omega_2\hat\delta(\omega_1)),
$$
for all $\omega_1,\omega_2\in \Omega^1$.
\endproclaim
\demo{proof}
We claim first that for every $a,x,y\in A$ we have
$$
\sigma_L(dy\,a\,dx) = \sigma_{L^*}(a\,dx\,d\delta(y)).\tag{6.12}
$$
Indeed, using Lemma 2.5 and the remarks following it, the
left side of (6.12) can be written
$$
\align
\sigma_L(dy\,a\,dx) &= \rho\circ\theta_L(dy\,a\,dx) =
\rho(L(yax)-yL(ax)-L(ya)x+yL(a)x) \\
&= -\rho(yL(ax)) - \rho(L(ya)x) + \rho(yL(a)x).
\endalign
$$
On the other hand,
$$
\sigma_{L^*}(a\,dx\,d\delta(y)) =
\rho(aL^*(x\delta(y))-axL^*(\delta(y))-aL^*(x)\delta(y)
+ axL^*(\bold 1)\delta(y)).
$$
Using $L^*(\bold 1) =0$ and the formula
$\rho(uL^*(v)) = \rho(L(u)v)$ the right side of the preceding
formula becomes
$$
\rho(L(a)x\delta(y))-\rho(L(ax)\delta(y))-\rho(yaL^*(x))
=\rho(yL(a)x) -\rho(yL(ax))-\rho(L(ya)x).
$$
It follows that
$$
\sigma_L(dy\,a\,dx) -\sigma_{L^*}(a\,dx\,d\delta(y)) = 0,
$$
as asserted in (6.12).

We may conclude that for all $\omega_2\in \Omega^1$ we have
$$
\sigma_L(dy\,\omega_2)= \sigma_{L^*}(\omega_2\hat\delta(dy)).
$$
Now for any $b,y\in A$ we have by Lemma 6.10
$$
\sigma_L(b\,dy\,\omega_2) = \sigma_L(dy\,\omega_2\,\delta(b)),
$$
and by (6.12) the right side is
$\sigma_{L^*}(\omega_2\,\delta(b)\,d\delta(y))$.  Thus for any
$\omega_1$ which is a finite sum of elements of the form
$b\, dy$ we have
$$
\sigma_L(\omega_1\omega_2) = \sigma_{L^*}(\omega_2\, \hat\delta(\omega_1)),
$$
and 6.11 follows\qed
\enddemo

Turning now to the proof of 6.9, we show
(ii) $\implies$ (i) $\implies$ (iii) $\implies$ (iv) $\implies$
(ii).
The first of these implications is immediate from Proposition 5.7.

\demo{(i) $\implies$ (iii)}
Proposition 5.7 implies that
$\omega_L(x,y) = \rho(xD(y))$ where $D$ is a derivation.  On the
other hand, by (6.11) we have
$$
\omega(x,y) = \rho(xL(y)) - \rho(L(x)y) = \rho(x(L(y)-L^*(y)).
$$
Hence $L-L^* = D$ is plainly a derivation.
\enddemo

\demo{(iii) $\implies$ (iv)}
Since $L-L^*$ is a first order differential operator its
symbol must vanish, and hence
$\sigma_L - \sigma_{L^*} = \sigma_{L-L^*} = 0$.  Thus by Lemma
6.11 we have
$$
\sigma_L(\omega_1\omega_2) =
\sigma_{L^*}(\omega_2\hat\delta(\omega_1)) =
\sigma_L(\omega_2\hat\delta(\omega_1)),
$$
as asserted.
\enddemo

\demo{(iv) $\implies$ (ii)}
Assuming the KMS condition 6.9 (iv),  the preceding argument can
clearly be reversed to show that $\sigma_{L-L^*} = 0$,
hence $D = L-L^*$ is a first order differential operator.  Since
$L(\bold 1) = L^*(\bold 1) = 0$ we have $D(\bold 1) = 0$ and thus
$D$ is a derivation.

Choose an operator $v\in A$ such that $D(x) = [v,x]$, $x\in A$.
We claim that $v+v^*$ is a scalar.  To see this we use the fact that
$D = L-L^*$ is a skew-adjoint operator on $L^2(A,\rho)$ which
satisfies $\rho\circ D = 0$.  Setting $\overline{D}(x) = D(x^*)^*$
we find that $\overline{D}$ and $D^*$ are related by
$$
\rho(\overline{D}(a)b) = \rho(aD^*(b)).
$$
Since $D^* = -D$ and since $\rho\circ D = 0$, we have
$$
\rho(aD^*(b)) = -\rho(aD(b)) = -\rho(D(ab)) +\rho(D(a)b) = \rho(D(a)b).
$$
Hence $\overline{D} = D$.  Using the fact that $D = [v,\cdot]$ and
$\overline{D} = [-v^*,\cdot]$
 we conclude that $[v+v^*,a] = 0$ for every $a\in A$,
hence $v+v^*$ is a scalar.

By replacing $v$ with $v-\lambda\bold 1$ for a suitable real number $\lambda$,
we can assume that $v^* = -v$; and the condition $\rho\circ[v,\cdot] = 0$
implies that $\delta(v) = v$, completing the proof of Theorem 6.9\qed
\enddemo

\proclaim{Corollary 6.13}
Let $(P,\<\cdot,\cdot\>)$ be a momentum space with Laplacian
$\Delta$ and let $L\in \Cal D(A,\rho)$ have the form
$$
L = \Delta + [v,\cdot]
$$
where $v^* = -v$, $\delta(v) = v$.  Then $L$ is exact.
\endproclaim
\demo{proof}
By Proposition 4.5, $L-L^* = 2[v,\cdot]$ is a derivation and thus
$L$ satisfies condition (iii) of Theorem 6.9\qed
\enddemo

\remark{Remarks}
We emphasize that the symbol $\sigma_L$ of an operator $L\in \Cal D(A,\rho)$
does not normally give rise to a modular cocycle.  More precisely,
let $\phi:A^3\to \Bbb C$ be the trilinear form
$$
\phi(a,x,y) = \sigma_L(a\,dx\,dy).
$$
A straightforward computation shows that
$$
b_\delta \phi = 0.
$$
However, $\phi$ does not satisfy the functional equation for
two-cochains
$$
\phi(a,x,y) = \phi(\delta^{-1}(y),a,x)
$$
except in rather special circumstances.  We omit
further discussion of the
latter since it does not come to bear on the sequel.

For an elliptic operator $L\in \Cal D(A,\rho)$, the form
$\omega_L(x,y) = \rho(xL(y)) - \rho(L(x)y)$ plays a role closely
analogous to the ``driving force" of a classical mechanical system.
Indeed, in what follows, the condition 6.9 (ii) for exactness,
$$
\omega_L(x,y) = -\rho(x[y,v])
$$
will occupy a position parallel to the hypothesis of classical mechanics
$$
F = -\text{grad } V
$$
that the driving force should be conservative.
\endremark

\example{Example 6.14}
We conclude this section by describing some elementary
examples of elliptic operators that are not exact.
Let $(A,\tau)$ be the pair consisting of a full
matrix algebra $M_n(\Bbb C)$, $n\geq 3$, with normalized trace
$\tau$.  Let $\alpha$ be a $*$-automorphism of $A$ and
consider the operator
$$
L(x) = \alpha(x) - x, \qquad x\in A.  \tag{6.15}
$$
One may verify directly that $L$ is the generator of the
semigroup
$$
\phi_t = e^{-t}\exp(t\alpha),\qquad t\geq 0
$$
and that $(\phi,\tau)$ is a Markov semigroup.  Thus,
$L$ is an elliptic operator in $\Cal D(A,\tau)$.
The following result implies that operators of
the form (6.15) are typically not exact.

\proclaim{Proposition}
Assume that there is an abelian $*$-subalgebra
$N\subseteq A$ such that $\alpha(N) = N$ and
$\alpha^2\restriction_N$ is not the identity map of $N$.
Then $L$ is not exact.
\endproclaim

\demo{proof}
Let $\omega_L:A\times A\to \Bbb C$ be the form
$$
\omega_L(x,y) = \tau(xL(y))-\tau(L(x)y).
$$
We will show that there is no derivation $D$ of $A$ for
which $\omega_L(x,y) = \tau(xD(y))$, $x,y\in A$.  This
means that $L$ fails to satisfy condition (ii) of Theorem
6.9, hence $L$ is not exact.

Using the fact that every automorphism of $A$ preserves
the trace, we have
$$
\align
\omega_L(x,y) &= \tau(x(\alpha(y)-y) -(\alpha(x)-x)y)
=\tau(x\alpha(y)-\alpha(x)y) \\
&= \tau(x(\alpha(y)-\alpha^{-1}(y)).
\endalign
$$
Thus we have to show that the operator $\alpha - \alpha^{-1}$ is
not a derivation.

But $\alpha-\alpha^{-1}$ leaves the
abelian $*$-subalgebra $N$ invariant,
so if it is a derivation it must be $0$ on N.  This
implies that $\alpha(x) = \alpha^{-1}(x)$ for every
$x\in N$, contradicting the hypothesis on $\alpha$\qed
\enddemo
\endexample

\subheading{7.  Classification of elliptic operators}

We now take up the classification of elliptic
operators in $\Cal D(A,\rho)$.  This is based on the classification
of {\it metrics}.  $A$ will denote a full matrix algebra.

\proclaim{Definition 7.1}
A metric on $A$ is a linear functional $g:\Omega^2\to \Bbb C$
satisfying
$$
g(\omega^*\omega) \geq 0,\qquad \omega\in \Omega^1.
$$
\endproclaim

A metric gives rise to a positive semidefinite inner product
$(\cdot,\cdot) : \Omega^1\times\Omega^1 \to \Bbb C$ by way of
$$
(\omega_1,\omega_2) = g(\omega_2^*\omega_1),
\qquad \omega_1,\omega_2\in \Omega^1
$$
and we have
$$
(a\omega_1,\omega_2) = (\omega_1,a^*\omega_2).  \tag{7.2}
$$
Conversely, any (semidefinite) inner product $(\cdot,\cdot)$
satisfying (7.2) arises in this way from a unique metric
$g:\Omega^2\to \Bbb C$.  We have found it more convenient to
work with metrics {\it qua} linear functionals rather than
with metrics {\it qua} inner products.  Notice, for example, that
the symbol of any elliptic operator $L$ defines a metric
$g$ by way of $g = -\sigma_L$.

Metrics can be concretely presented as follows.  Let
$z_1,z_2,\dots,z_r \in A$ and set
$$
g(a\, dx\,dy) = \sum_{k=1}^r\rho(a[z_k^*,x][z_k,y]).  \tag{7.3}
$$
It is easily seen that $g$ is a metric.  Notice that this metric
satisfies
$$
g(a\xi) = g(\xi\,\delta(a)),\qquad a\in A,\xi \in \Omega^2,  \tag{7.4}
$$
$\delta$ denoting the modular automorphism of $\rho$.  Conversely, it
can be shown that any metric satisfying (7.4) can be expressed in the
form (7.3) for some set of elements $z_1,z_2,\dots,z_r$ in $A$.
Note too that the $z_k$ do not need to be
self-adjoint or skew-adjoint, and they need not satisfy
$\delta(z_k) = z_k$.

In this section we consider metrics $g$ which satisfy the KMS condition
of Theorem 6.9 (iv),
$$
g(\omega_1\omega_2)  = g(\omega_2\, \hat\delta(\omega_1)),
\qquad \omega_1,\omega_2\in \Omega^1.   \tag{7.5}
$$
We prove that a KMS metric can be
decomposed into a sum of the form (7.3) in which the operators
$z_1,\dots,z_r$ satisfy the additional conditions
$z_k^* = z_k$ and $\delta(z_k) = z_k$ for all $k$, and
that moreover there is a unique momentum space associated
with $g$.

\proclaim{Theorem 7.6}
Let $g$ be a nonzero metric which satisfies the KMS condition
\rm{(7.5)}.
There is a linearly independent set of self-adjoint
operators $x_1,\dots,x_n\in A$ such that $\delta(x_k) = x_k$
and $\rho(x_k) = 0$ for every $k$, and for which
$$
g(a\,dx\,dy) = \sum_{k=1}^n\rho(a[x_k,x][x_k,y])
$$
for every $a,x,y\in A$.

If $x_1^\prime,\dots,x_m^\prime$ is another finite set satisfying
all of these conditions then $m=n$ and there is a real orthogonal
$n\times n$ matrix $(u_{ij})$ such that
$$
x_k^\prime = \sum_{j=1}^nu_{kj}x_j,\qquad 1\leq k\leq n.
$$
\endproclaim

\remark{Remarks}
Before giving the proof of Theorem 7.6 we discuss some
immediate consequences.

Notice that by replacing $x_k$ with $p_k = \sqrt{-1}\,x_k$ we
obtain a linearly independent
set of skew-adjoint elements $\{p_1,\dots,p_n\}$
satisfying $\rho([p_k,a]) = \rho(p_k) =0$
for every $a\in A$ such that
$$
g(a\,dx\,dy) = -\sum_{k=1}^n\rho(a[p_k,x][p_k,y]),
$$
and for which a similar uniqueness holds.  Thus the real
vector space $P$ spanned by $\{p_1,\dots,p_n\}$ is independent
of the particular choice of $p_1,\dots,p_n$, as is the
inner product defined on $P$ by
$$
\<p_i,p_j\> = \delta_{ij}.
$$
We conclude that {\it every KMS metric $g$ is associated with
a unique momentum space $(P,\<\cdot,\cdot\>)$}.

This association $g\rightsquigarrow (P,\<\cdot,\cdot\>)$ of a
momentum space with a KMS metric is in fact bijective.  Indeed,
if we let $\Delta$ be the Laplacian of $(P,\<\cdot,\cdot\>)$ then
from (2.10) we find that $2g=-\sigma_\Delta$ is determined by
$\Delta$, and therefore by $(P,\<\cdot,\cdot\>)$.
Thus Theorem 7.6 implies the following
\endremark

\proclaim{Corollary 1}
The above association defines a bijective correspondence between
the set of KMS metrics and the set of momentum spaces in $(A,\rho)$.
\endproclaim

Indeed, Theorem 7.6 can be restated in the following equivalent
(and somewhat more invariant) way.

\proclaim{Theorem 7.6A}
For every KMS metric
$g:\Omega^2\to \Bbb C$ there is a unique
momentum space $(P,\<\cdot,\cdot\>)$ such that $2g = -\sigma_\Delta$,
$\Delta$ denoting the Laplacian of $(P,\<\cdot,\cdot\>)$.
\endproclaim

Two metrics $g_1, g_2$ are said to be {\it equivalent} if there is a
$*$-automorphism $\alpha:A\to A$ such that $g_2(\xi) = g_1(\hat\alpha(\xi)),
\xi\in \Omega^2$, $\hat\alpha$ denoting the induced mapping of
$\Omega^2$:
$$
\hat\alpha(a\,dx\,dy) = \alpha(a)\,d\alpha(x)\,d\alpha(y).
$$
The above remarks lead immediately to the following classification
of KMS metrics in terms of their momentum spaces.

\proclaim{Corollary 2}
Let $g_1,g_2$ be two KMS metrics with respective momentum spaces
$(P_1,\<\cdot,\cdot\>_1)$ and $(P_2,\<\cdot,\cdot\>_2)$.  $g_1$ and
$g_2$ are equivalent iff there is a unitary operator $u\in A$
such that
$$
\align
uP_1u^* &= P_2,\qquad \text{and}\\
\<upu^*,uqu^*\>_2 &= \<p,q\>_1,\qquad p,q\in P_1.
\endalign
$$
\endproclaim
\demo{proof}
Assuming that $g_1$ and $g_2$ are equivalent, let $\alpha$ be a
$*$-automorphism of $A$ such that $g_2=g_1\circ\hat\alpha$.
We may find a unitary operator $u\in A$ such that $\alpha(x) = uxu^*$.
We can obviously define a momentum space
$(uP_1u^*,\<\cdot,\cdot\>^\prime)$ where
$\<upu^*,uqu^*\>^\prime=\<p,q\>_1$,
$p,q\in P_1$; and
$(uP_1u^*,\<\cdot,\cdot\>^\prime)$ gives rise to the metric
$g_1\circ\hat\alpha = g_2$.  Hence by the uniqueness assertion
of Theorem 7.6 we must have
$(P_2,\<\cdot,\cdot\>) = (uP_1u^*,\<\cdot,\cdot\>^\prime)$.

The opposite implication is equally straightforward\qed
\enddemo

\demo{proof of Theorem 7.6}
Since $g(\omega^*\omega) \geq 0$ for every $\omega\in \Omega^1$,
a standard construction leads to a finite dimensional Hilbert
space $H_g$ and a complex linear map $\theta:\Omega^1\to H_g$
satisfying

$$
\align
\theta_g(\Omega^1) &= H_g,\qquad \text{and} \tag{7.7.1}\\
g(\omega_2^*\omega_1) &= \<\theta(\omega_1),\theta(\omega_2)\>,
\qquad \omega_1,\omega_2\in \Omega^1.  \tag{7.7.2}
\endalign
$$
Let $a$ be an element of $A$.
The formula $g(\omega_2^*a\omega_1) = g((a^*\omega_2)^*\omega_1)$
implies that
$$
\<\theta(a\omega_1),\theta(\omega_2)\> =
\<\theta(\omega_1),\theta(a^*\omega_2)\>\qquad \omega_i\in \Omega^1.
$$
It follows that there is a unique $*$-representation
$\pi:A\to \Cal B(H_g)$  defined by
$$
\pi(a)\theta(\omega) = \theta(a\omega), \qquad \omega\in \Omega^1.  \tag{7.8}
$$
$\pi(\bold 1) = \bold 1$ because $\bold 1\cdot\omega = \omega$ for every
$\omega\in \Omega^1$.

Next, we show that there is a natural conjugation $J$ of
$H_g$, that is, an antilinear operator on $H_g$ satisfying
$J^2 = \bold 1$ and $\<J\xi,J\eta\> = \<\eta,\xi\>$, for all
$\xi,\eta\in H_g$.  For that we require some information about
the behavior of $g$ with respect to the (complex) modular
automorphism group of the distinguished state on $A$,
$\rho(a) = \text{trace}(ha)$.  For $z\in \Bbb C$ there is an
automorphism $\sigma_z$ of the algebra structure of $A$ defined
by
$$
\sigma_z(a) = h^{-iz}ah^{iz}.
$$
We have $\sigma_z\sigma_w= \sigma_{z+w}$ and
$\sigma_z(a^*) = \sigma_{\bar z}(a)^*$. For fixed
$a\in A$, $z\mapsto \sigma_z(a)$ defines an entire
function from $\Bbb C$ to the Banach space $A$ which
is uniformly bounded on horizontal strips
$-M\leq \Im(z)\leq M$, $0<M<\infty$.  For real $z$,
$\sigma_z$ is a $*$-automorphism of $A$.

 Each automorphism
$\sigma_z$ determines a natural automorphism $\hat{\sigma_z}$
of the differential algebra $\Omega^*$, and in particular
$\hat{\sigma_z}$ acts on $\Omega^1$ and $\Omega^2$ by
$$
\align
\hat\sigma_z(a\,dx) &= \sigma_z(a)\,d\sigma_z(x),\\
\hat\sigma_z(a\,dx\,dy) &= \sigma_z(a)\,d\sigma_z(x)\,d\sigma_z(y).
\endalign
$$
If $z=t$ is real then
$\hat{\sigma_t}(\omega^*) = \hat{\sigma_t}(\omega)^*$
for all $\omega$.
The modular automorphism $\delta$ introduced in \S 5 is
given by $\delta(a) = \sigma_i(a), a\in A$.

We now define a second involution $\#$ of $\Omega^1$ using the
natural square root of the modular automorphism
$$
\delta^{1/2}(a) = \sigma_{i/2}(a) = h^{1/2}ah^{-1/2}
$$
as follows:
$$
(a\,dx)^\# = \hat{\delta^{1/2}}((a\,dx)^*) =
-d\delta^{1/2}(x^*)\,\delta^{1/2}(a^*).  \tag{7.9}
$$
Since
$\delta^{1/2}(x^*) = \delta^{-1/2}(x)^*$ we see
that $\omega\mapsto\omega^\#$ defines an involution
in $\Omega^1$.

\proclaim{Lemma 7.10}
Every KMS metric $g:\Omega^2\to \Bbb C$ satisfies
$$
\align
g(\hat{\sigma_z}(\xi)) &= g(\xi), \qquad \xi\in \Omega^1, z\in \Bbb C
\tag{i}\\
g(\omega_2^*\omega_1) &= g((\omega_1^\#)^*\omega_2^\#)
\qquad \omega_1,\omega_2\in \Omega^1,\tag{ii}
\endalign
$$
We have $(a\omega)^\# = \omega^\#\delta^{1/2}(a^*)$
and $(\omega b)^\# = \delta^{1/2}(b^*)\omega^\#$,
for every $a,b\in A$, $\omega\in \Omega^1$.
\endproclaim

\demo{proof}

For (i), notice first that $g\circ\hat\delta = g$.  Indeed,
for $\omega_1,\omega_2\in \Omega^1$ we can apply the KMS
condition twice to obtain
$$
g(\hat\delta(\omega_1\omega_2)) =
g(\hat\delta(\omega_1)\hat\delta(\omega_2)) =
g(\omega_2\hat\delta(\omega_1)) = g(\omega_1\omega_2),
$$
and $g\circ\hat\delta = g$ follows because $\Omega^2$ is
spanned by products of the form $\omega_1\omega_2$, $\omega_k\in \Omega^1$.

Fix an element $\xi\in \Omega^2$ and consider the entire function
$$
f(z) = g(\hat\sigma_z(\xi)),\qquad z\in \Bbb C.
$$
By the preceding remarks $f$ is entire and bounded on the horizontal
strip $0\leq \Im(z)\leq 1$.  We claim that $f(z+i) = f(z)$ for every
$z\in \Bbb C$.  Indeed, since
$\sigma_{z+i} = \sigma_i\sigma_z = \delta\sigma_z$ and since
$g\circ\hat\delta = g$ we have
$$
f(z+i) = g(\hat\delta(\hat\sigma_z(\xi))) = g(\hat\sigma_z(\xi)) = f(z).
$$
Thus $f$ is a bounded entire function which, by Liouville's theorem, must
be a constant.

For (ii) we write
$$
g((\omega_1^\#)^*\omega_2^\#) =
g((\hat{\delta^{1/2}}(\omega_1^*))^*\hat{\delta^{1/2}}(\omega_2^*))
= g(\hat{\delta^{-1/2}}(\omega_1)\hat{\delta^{1/2}}(\omega_2^*)).
$$
By what was just proved we can use
$g\circ\hat\delta^{1/2} = g$ on the right hand term to obtain
$$
g(\omega_1\delta(\omega_2^*)) = g(\omega_2^*\omega_1),
$$
and (ii) follows.  The formula $(a\omega)^\# = \omega^\#\delta^{1/2}(a^*)$
follows directly from the definition of $^\#$.  \qed
\enddemo

Lemma 7.10 (ii) implies that
$$
\<\theta(\omega_1),\theta(\omega_2)\> =
\<\theta(\omega_2^\#),\theta(\omega_1^\#)\>,
$$
and hence we can define a unique antilinear isometry $J$ of $H_g$
by
$$
J\theta(\omega) = \theta(\omega^\#),\qquad \omega \in \Omega^1. \tag{7.11}
$$
Since $\omega^{\#\#} = \omega$ it follows that $J^2=\bold 1$, and
hence $J$ is a conjugation of $H_g$.

We now define a unital $*$-antirepresentation
$\overset\circ\to\pi:A\to \Cal B(H_g)$ by
$$
\overset\circ\to \pi(a) = J\pi(a^*)J.
$$
{}From the definitions of $\pi$ and $J$ we find that

$$
\overset\circ\to \pi(a)\theta(\omega) = \theta(\omega\,\delta^{1/2}(a)),
\qquad a\in A,\omega\in \Omega^1,\tag{7.12}
$$
and hence $\pi(a)\overset\circ\to\pi(b) = \overset\circ\to\pi(b)\pi(a)$
for every $a,b\in A$.  Thus we have a $*$-representation of the
\cstar\ $A\otimes A^o$ ($A^o$ denoting the \cstar\ opposite to $A$)
by way of
$$
\pi\otimes\overset\circ\to \pi:a\otimes b\in A\otimes A^o \mapsto
\pi(a)\overset\circ\to \pi(b).
$$

Finally, the action of the modular group
$\hat{\sigma_t}:\Omega^1\to \Omega^1$ can be implemented by a one-parameter
group of unitary operators $U = \{U_t: t\in \Bbb R\}$,
which is defined on $H_g$ as follows:
$$
U_t\theta(\omega) = \theta(\hat{\sigma_t}(\omega)),
\qquad t\in \Bbb R,\omega\in \Omega^1.
$$
Indeed, for $\omega_1,\omega_2\in \Omega^1$,
$$
\align
\<\theta(\hat{\sigma_t}(\omega_1)),\theta(\hat{\sigma_t}(\omega_2)\>
&= g((\hat{\sigma_t}(\omega_2))^*\hat{\sigma_t}(\omega_1) )
= g((\hat{\sigma_t}(\omega_2^*)\hat{\sigma_t}(\omega_1)) \\
&= g(\hat{\sigma_t}(\omega_2^*\omega_1))
= g(\omega_2^*\omega_1) = \<\theta(\omega_1),\theta(\omega_2)\>,
\endalign
$$
and thus each $U_t$ is well-defined by the above formula.
$U$ obeys the group property and is strongly continuous.

\proclaim{Lemma 7.13}
We may assume that $H_g$ is the Hilbert space $A^n$ consisting of all
$n$-tuples $\bar z = (z_1,\dots,z_n)$, $z_k\in A$, with inner
product
$$
\<\bar z,\bar w\> = \sum_{k=1}^n\rho(w_k^*z_k),
$$
and that there is a set of $n$ real numbers
$\lambda_1,\dots,\lambda_n\in \Bbb R$ such that the quadruple
$\pi$, $\overset\circ\to\pi$, $U$, $J$ acts in the following way
$$
\align
\pi(a)\bar z &= (az_1,\dots,az_n) \tag{7.13.1}\\
\overset\circ\to \pi(b)\bar z &=
(z_1\delta^{1/2}(b),\dots,z_n\delta^{1/2}(b)). \tag{7.13.2}\\
U_t\bar z &= (e^{i\lambda_1 t}\sigma_t(z_1), \dots, e^{i\lambda_n
t}\sigma_t(z_n)).
\tag{7.13.3}\\
J\bar z &= (\delta^{1/2}(z_1^*),\dots,\delta^{1/2}(z_n^*)), \tag{7.13.4}\\
\endalign
$$
\endproclaim
\remark{Remark}
We point out that in the course of the argument below, we will
eventually prove that $\lambda_1 = \dots = \lambda_n = 0$.
\endremark
\demo{proof}
Since $A\otimes A^o$ is another matrix algebra the commutant of the
set of operators $\pi(A)\cup \overset\circ\to \pi(A^o)$ is a type
$I_n$ factor for some $n=1,2,\dots$, and the representation
$\pi\otimes\overset\circ\to \pi$ is unitarily equivalent to a direct
sum of $n$ copies of any irreducible representation of $A\otimes A^o$.

The standard representation $\lambda:A\to \Cal B(L^2(A,\rho))$ is
defined by
$$
\lambda(a)z=az, \qquad z\in A
$$
and there is a conjugation $J_0$ of $L^2(A,\rho)$ defined by
$$
J_0(z) = \delta^{1/2}(z^*).
$$
One has $J_0\lambda(A)J_0 = \lambda(A)^\prime$.  The corresponding
antirepresentation is defined by
$$
\overset\circ\to \lambda(b) = J_0\lambda(b^*)J_0,
$$
and one computes that
$$
\overset\circ\to \lambda(b)z = z\delta^{1/2}(b).
$$
Thus $\lambda\otimes \overset\circ\to \lambda$ defines an irreducible
representation of $A\otimes A^o$ on the Hilbert space $L^2(A,\rho)$.
Hence $\pi\otimes \overset\circ\to\pi$ is unitarily equivalent to
a direct sum of $n$ copies of $\lambda\otimes\overset\circ\to\lambda$.

After this change in coordinates, we may assume that $\pi$ and
$\overset\circ\to\pi$ act on $A^n$ as specified in (7.13.1) and
(7.13.2).

We show next that we can also achieve (7.13.3).  Let
$\overset\circ\to U_t$ be the unitary group which acts on
$A^n$ by way of
$$
\overset\circ\to U_t\bar z = (\sigma_t(z_1),\dots,\sigma_t(z_n)),
\qquad t\in \Bbb R, z_i\in A.
$$
Since both $U_t$ and $\overset\circ\to U_t$ induce the same action
on $\pi(A)$ and $\overset\circ\to \pi(A)$, namely
$$
\align
U_t\pi(a)U_t^* &= \overset\circ\to U_t\pi(a)\overset\circ\to U_t^*
=\pi(\sigma_t(a)),\\
U_t\overset\circ\to\pi(a)U_t^* &=
\overset\circ\to U_t\overset\circ\to\pi(a)\overset\circ\to U_t^*
=\overset\circ\to \pi(\sigma_t(a))
\endalign
$$
it follows that $W_t = U_t\overset\circ\to U_t^*$ commutes with
all operators in $\pi(A)\cup \overset\circ\to \pi(A)$, and we have
$U_t = W_t\overset\circ\to U_t$.

Notice that $W = \{W_t: t\in \Bbb R\}$ is a one-parameter unitary group
in the commutant of $\pi(A)\cup \overset\circ\to \pi(A)$.  Indeed,
since $\sigma_t(a) = h^{-it}ah^{it}$ is an inner automorphism of
$A$ for each $t$ it follows that $\overset\circ\to U_t$ belongs
to the von Neumann algebra generated by $\pi(A)\cup \overset\circ\to \pi(A)$.
Thus
$$
W_{s+t}\overset\circ\to U_{s+t} = U_{s+t} = U_sU_t =
W_s\overset\circ\to U_sW_t\overset\circ\to U_t =
W_sW_t\overset\circ\to U_s\overset\circ\to U_t =
W_sW_t\overset\circ\to U_{s+t},
$$
so that $W_sW_t = W_{s+t}$.  By the spectral theorem
there are mutually orthogonal minimal projections
$E_1,\dots,E_n$ in $(\pi(A)\cup \overset\circ\to\pi(A))^\prime$
and real numbers $\lambda_1,\dots,\lambda_n$ such that
$$
W_t = \sum_{k=1}^n e^{i\lambda_kt}E_k.
$$
Thus there is a unitary operator $V$ in the commutant of
$\pi(A)\cup\overset\circ\to \pi(A)$ which brings the unitary
group $W = \{W_t: t\in \Bbb R\}$ into diagonal form relative
to the coordinates given in (7.13.3).  Conjugation by this operator
$V$ does not change $\pi$ or $\overset\circ\to \pi$, and thus we have
achieved (7.13.1), (7.13.2) and (7.13.3) simultaneously.

Let $J_n$ be the natural conjugation of $A^n$, defined by
$$
J_n(\bar z) = (\delta^{1/2}(z_1^*),\dots,\delta^{1/2}(z_n^*)),
$$
and let $J$ be the conjugation defined by (7.11).  It remains to
 show that there
is a unitary operator $W$ on $A^n$ such that
$$
\align
W &\in (\pi(A) \cup \overset\circ\to\pi(A)\cup \{U_t: t\in \Bbb R\})^\prime,
\qquad \text{and} \tag{7.14}\\
J &= WJ_nW^{-1}.\tag{7.15}
\endalign
$$
Once we have such a $W$, we may use it to bring $J$ into the
form (7.13.4) without disturbing what has already been achived
in (7.13.1), (7.13.2) and (7.13.3).

$W$ is defined as follows.  The standard conjugation $J_n$ of $A^n$ satisfies
$$
\align
J_n(\pi(a^*)\bar z) &= J_n(a^*z_1,\dots,a^*z_n) =
(\delta^{1/2}(z_1^*a),\dots,\delta^{1/2}(z_n^*a)) \\
&=(\delta^{1/2}(z_1)\delta^{1/2}(^*a),\dots,\delta^{1/2}(z_n)\delta^{1/2}(^*a))
=\overset\circ\to\pi(a)J_n(\bar z),
\endalign
$$
and hence
$$
J_n\pi(a^*) J_n = \overset\circ\to\pi(a).  \tag{7.16}
$$
Since $J_n^2 = \bold 1$ we also have
$$
J_n\overset\circ\to\pi(a)J_n = \pi(a^*).  \tag{7.17}
$$
Since $J$ also satisfies (7.16) and (7.17) it follows that the
unitary operator $V = JJ_n$ must commute with both operator algebras
$\pi(A)$ and $\overset\circ\to\pi(A)$, and we have $J=VJ_n$.

Notice
that $V$ also commutes with the unitary group $\{U_t: t\in \Bbb R\}$.
Indeed, from the definitions of $U_t$ and $J$ we have
$$
JU_t\theta(\omega) = J\theta(\hat{\sigma_t}\omega)
= \theta((\hat{\sigma_t}\omega)^\#).
$$
But for real $t$ we have
$$
(\hat{\sigma_t}\omega)^\# = \hat\sigma_{i/2}((\hat\sigma_t\omega)^*)
= \hat\sigma_{i/2}\hat\sigma_t(\omega^*) =
\hat\sigma_t(\hat\sigma_{i/2}(\omega^*)) = \hat\sigma_t(\omega^\#).
$$
Hence the right side of the previous equation is
$$
\theta(\hat\sigma_t(\omega^\#)) = U_t\theta(\omega^\#)= U_tJ\theta(\omega).
$$
We conclude that $JU_t = U_tJ$.  On the other hand,
from the representation of $U_t$ achieved in (7.13.3) we find that
$$
\align
J_nU_t\bar z &=
J_n(e^{i\lambda_1t}\sigma_t(z_1),\dots,e^{i\lambda_nt}\sigma_t(z_n))\\
&= (\sigma_{i/2}(e^{-i\lambda_1t}\sigma_t(z_1^*),
\dots,\sigma_{i/2}(e^{-i\lambda_nt}\sigma_t(z_n^*) \\
&= (\sigma_t((e^{i\lambda_1t}\sigma_t(z_1))^\#),\dots,
\sigma_t((e^{i\lambda_nt}\sigma_t(z_n))^\#))
=U_tJ_n\bar z.
\endalign
$$
Thus $U_t$ commutes with both $J$ and $J_n$, and therefore it
must commute with $V = JJ_n$.

We now show that any ``reasonable" square root $W$ of $V$ has the
properties (7.14) and (7.15).  Indeed, let
$$
V = \sum_{k=1}^r\lambda_kE_k
$$
be the spectral decomposition of $V$, where $\lambda_1,\dots,\lambda_r$ are
the distinct eigenvalues of $V$ and $E_1,\dots,E_r$ are its minimal spectral
projections.  Set
$$
W = \sum_{k=1}^r\lambda^{1/2}_kE_k,
$$
where $\lambda^{1/2}_k$ denotes either square root of $\lambda_k$.
It is clear that $W$ commutes with both $\pi(A)$ and $\overset\circ\to\pi(A)$,
and $W$ commutes with each $U_t$ because it belongs to the
\cstar\ generated by the operator $V\in \{U_t: t\in \Bbb R\}^\prime$.
Since $(UJ_n)^2 = J^2 = \bold 1$ it follows that $J_nUJ_n = U^{-1}$.
Now the operator mapping
$$
\alpha:T \mapsto J_n T^* J_n
$$
is a $*$-antiautomorphism of the von Neumann algebra $\Cal B(A^n)$ which
carries the abelian \cstar\ $C^*(U)$ generated by $U$ onto itself in such
a way that $\alpha(U) = U$.  It follows that $\alpha(T) = T$ for all
$T\in C^*(U)$.  In particular $\alpha(W) = W$, i.e., $J_nW^{-1}J_n = W$,
and therefore $U = W^2 = WJ_nW^{-1}J_n$.  Finally, we see that
$J = UJ_n = WJ_nW^{-1}$, as required\qed
\enddemo

We now describe $\theta$ in terms of the coordinates of
Lemma 7.13.

\proclaim{Lemma 7.18}
Considering $A^n$ as a bimodule over $A$ in the usual way,
$$
a\bar z b = (az_1b,\dots,az_nb),\qquad \bar z\in A^n,
$$
we have $\theta(a\omega b) = a\theta(\omega)b$, $a,b\in A$,
$\omega\in \Omega^1$.
\endproclaim
\demo{proof}
We have $\theta(a\omega) = \pi(a)\theta(\omega) = a\theta(\omega)$,
by (7.13.1).  For the right action we have
$$
\align
\theta(\omega b) &= J\theta((\omega b)^\#) =
J\theta(\delta^{1/2}(b^*)\hat{\delta^{1/2}}(\omega^*))
=J\pi(\delta^{1/2}(b^*))\theta(\omega^\#) \\
&=J\pi(\delta^{-1/2}(b)^*)J\theta(\omega)=
\overset\circ\to\pi(\delta^{-1/2}(b))\theta(\omega).
\endalign
$$
Notice that $\overset\circ\to\pi(\delta^{-1/2}(b))$ gives the
(untwisted) right module action of $b$, since by (7.13.2)
$$
\overset\circ\to \pi(\delta^{-1/2}(b))\bar z = (z_1b,\dots,z_nb),
$$
$\bar z\in A^n$.  Hence $\theta(\omega b) = \theta(\omega)b$\qed
\enddemo

It follows immediately from Lemma 7.18 that $x\mapsto \theta(dx)$ is
a derivation of $A$ into the bimodule $A^n$, and hence there are derivations
$D_1,\dots,D_n$ of $A$ such that
$$
\theta(a\,dx) = (aD_1(x),\dots,aD_n(x)),\qquad a,x\in A.
$$

\proclaim{Lemma 7.19}
The derivations $D_1,\dots,D_n$ are linearly independent.
\endproclaim
\demo{proof}
Suppose that $\lambda_1,\dots,\lambda_n \in \Bbb C$ are such that
$$
\sum_{k=1}^n\lambda_k D_k(x) = 0,\qquad x\in A.
$$
Then for every $a,x\in A$ we have
$$
\sum_{k=1}^n\lambda_k a D_k(x) = 0.
$$
It follows that the linear operator
$L:A^n\to A$ defined by $L(\bar z) = \sum _k\lambda_kz_k$ satisfies
$L(\theta(a\,dx)) = 0$ for all $a,x\in A$; hence $L=0$ on the range
of $\theta$.  Since $\theta(\Omega^1) = H_g = A^n$, $L$ must be $0$\qed
\enddemo

The remainder of the proof is devoted to showing that $D_k(a)^* = -D_k(a^*)$
and $\rho\circ D_k = 0$, $1\leq k\leq n$.  This is done in two steps, the
first of which is the following.

\proclaim{Lemma 7.20}The derivations $D_1,\dots,D_n$ have the
form $D_k(a) = [x_kh^{-1/2},a]$, where $x_k=x_k^*$ and $h$ is the
density matrix of the state $\rho$.
\endproclaim
\demo{proof}
Applying (7.13.3) to the $n$-tuple $\bar z= (D_1(x),\dots,D_n(x))$
and using $J\theta(dx) = \theta((dx)^\#) = -\theta(d\delta^{1/2}(x^*))$
we obtain
$$
-D_k(\delta^{1/2}(x^*)) = \delta^{1/2}(D_k(x)^*),\qquad 1\leq k\leq n.
$$
After choosing elements $z_k\in A$ so that $D_k(x) = [z_k,x]$ and
unravelling the preceding formula, one finds that for all
$x\in A$ one has
$$
-[z_k,\delta^{1/2}(x^*)] = -[\delta^{1/2}(z_k^*),\delta^{1/2}(x^*)].
$$
Thus $z_k-\delta^{1/2}(z_k^*)$ commutes with everything in $A$, and
must be a scalar multiple of $\bold 1$ for every $1\leq k\leq n$.

Now if $z\in A$ and $\lambda\in \Bbb C$ are such that
$z-\delta^{1/2}(z^*) = \lambda\bold 1$, then by applying the normalized
trace $\tau$ on $A$ we find that $\lambda = \tau(z) - \overline{\tau(z)}$
is imaginary.  Writing $\lambda = i\lambda_0$ where $\lambda_0\in \Bbb R$
we can replace $z$ with $z_0 = z - \frac{i\lambda_0}{2}\bold 1$ and we find
that $z_0 - \delta^{1/2}(z_0^*) = 0$.  Noting that
$\delta^{1/2}(x) = h^{1/2}xh^{-1/2}$, the latter implies that
$x=z_0h^{-1/2}$ is self-adjoint, and hence
$$
z = xh^{-1/2} + \frac{i\lambda_0}{2}\bold 1.
$$
The required representation of $z_1,\dots,z_n$ follows\qed
\enddemo

Now by definition of $U_t$ we have
$$
U_t(D_1(x),\dots,D_n(x)) = U_t(\theta(dx)) = \theta(d\sigma_t(x))
=(D_1(\sigma_t(x)),\dots,D_n(\sigma_t(x))),
$$
while from (7.13.3) we have
$$
U_t(D_1(x),\dots,D_n(x)) =
(e^{i\lambda_1t}\sigma_t(D_1(x)),\dots,e^{i\lambda_nt}\sigma_t(D_n(x)).
$$
It follows that
$$
\sigma_t(D_k(x)) = e^{-i\lambda_kt}D_k(\sigma_t(x)),
\qquad 1\leq k\leq n,t\in \Bbb R.
$$
By (7.20), $D_k$ has the form $D_k(x) = [x_kh^{-1/2},x]$, and thus
$$
\sigma_t(x_kh^{-1/2}) - e^{-i\lambda_kt}x_kh^{-1/2}
$$
must be a scalar.  The following lemma shows that this implies that
$x_k$ must commute with $h$.

\proclaim{Lemma 7.21}
Let $\lambda\in \Bbb R$ and let $x$ be a self adjoint element of
$A$ such that
$$
\sigma_t(xh^{-1/2}) - e^{i\lambda t}xh^{-1/2}
$$
is a scalar for every $t\in \Bbb R$.  Then $\lambda = 0$ and $x$
commutes with $h$.
\endproclaim
\demo{proof}
For every $t\in \Bbb R$ let $\mu(t)$ be the complex number defined
by
$$
\sigma_t(xh^{-1/2}) - e^{i\lambda t}x_kh^{-1/2} = \mu(t)\bold 1.
$$
Applying the normalized trace $\tau$ of $A$ to both sides of this
equation we find that for $c=\tau(xh^{-1/2})$ we have
$\mu(t) = c(1-e^{i\lambda t})$.  Thus
$$
\sigma_t(xh^{-1/2} - c\bold 1) = e^{i\lambda t}(xh^{-1/2} - c\bold 1).
$$
After multiplying through on the right by $h^{1/2}$ we find that
$$
\sigma_t(x - ch^{1/2}) = e^{i\lambda t}(x - ch^{1/2}), \qquad t\in \Bbb
R.\tag{7.22}
$$
Now $c = \tau(xh^{-1/2}) = \tau(h^{-1/4}xh^{-1/4})$ is a real
number because $x$ is a self-adjoint operator, and thus
the operator $x-ch^{1/2}$ appearing in (7.22) is self-adjoint.
The left side of (7.22) is therefore a self-adjoint operator for
every $t\in \Bbb R$, and hence (7.22) implies that $\lambda = 0$
except in the trivial case where $x = ch^{1/2}$.

In either case it follows that $x-ch^{1/2}$ is fixed under the action
of $\sigma_t$, $t\in \Bbb R$, and hence $x$ must commute with $h$  \qed
\enddemo

Applying Lemma 7.21 to the derivations $D_k$ we find that
$D_k(a) = [x_kh^{-1/2},a] = [h^{-1/4}x_kh^{-1/4},a]$ is implemented
by the self-adjoint operator $y_k = h^{-1/4}x_kh^{-1/4}$,
where $y_k$ commutes with $h$.  By replacing $y_k$
with $y_k - \rho(y_k)\bold 1$ we can also assume that
$\rho(y_k) = 0$.  That completes the proof of all but the
uniqueness assertion of Theorem 7.6.

For uniqueness, let $x_1^\prime,\dots,x_n^\prime$ be a linearly
independent set of self-adjoint operators satisfying
$\delta(x_k^\prime) = x_k^\prime$, $\rho(x_k^\prime) = 0$,
and
$$
g(a\,dx\,dy) = \sum_{k=1}^m \rho(a[x_k^\prime,x][x_k^\prime,y])\tag{7.23}
$$
for all $a,x,y\in A$.  Define a linear map
$\theta^\prime:\Omega^1 \to A^m$ by
$$
\theta^\prime(a\,dx) = (a[x_1^\prime,x],\dots,a[x_m^\prime,x]).  \tag{7.24}
$$
If we make $A^m$ into a Hilbert space $m\cdot L^2(A,\rho)$ as in
Lemma 7.13, and into a bimodule as in Lemma 7.18, then we find that
$$
\theta^\prime(a\omega b) = a\theta^\prime(\omega)b,
\qquad a,b\in A,\omega\in \Omega^1. \tag{7.25}
$$

We claim that $\theta^\prime(\Omega^1) = A^m$.  For that,
note that (7.25) implies that the subspace $\theta^\prime(\Omega^1)$
is invariant under left and right multiplication by elements of $A$,
and hence the projection $P$ onto the range of $\theta^\prime$ must
have the form $P = (p_{ij}\bold 1)$ where $(p_{ij})$ is an $m\times m$
matrix of complex scalars.  If $P\neq \bold 1$ then $\bold 1-P$ must
have a nonzero row with entries $\mu_1,\dots,\mu_m$.  Since
$(\bold 1 - P)\theta^\prime(dx) = 0$ for every $x\in A$ it follows that
$$
\mu_1[x_k^\prime,x] + \dots + \mu_m[x_m^\prime,x] = 0
$$
for every $x\in A$ and hence $\sum\mu_kx_k^\prime$ must be a
scalar multiple $\nu\bold 1$ of $\bold 1$.
However, since $\rho(x_k^\prime) = 0$ and $\rho$ is a state,
we must have $\nu = 0$.  Hence $\sum\mu_kx_k^\prime = 0$,
contradicting linear independence.

Now the formulas (7.23) and (7.24) imply that for all
$\omega_1$, $\omega_2\in \Omega^1$,
$$
\<\theta^\prime(\omega_1),\theta^\prime(\omega_2)\>
= g(\omega_2^\# \omega_1) = \<\theta(\omega_1),\theta(\omega_2))\>.
$$
Thus we can define a unique unitary operator $W:A^n\to A^m$ by
$$
W\theta(\omega) = \theta^\prime(\omega),\qquad \omega\in \Omega^1.
$$
Since $W$ is unitary and both $\theta$ and $\theta^\prime$ are
bimodule homomorphisms it follows that $W$ must implements
a $*$ isomorphism of the commutant of $\pi(A) \cup \overset\circ\to \pi(A)$
onto that of $\pi^\prime(A) \cup \overset\circ\to \pi^\prime(A)$. Since
these commutants are factors of type $I_n$ and $I_m$ respectively,
we conclude that $m=n$, and that $W$ commutes with both the left
and right module actions of $A$ on $A^n$.  It follows that there is
an $n\times n$ scalar unitary matrix $(w_{ij})$ such that
$W = (w_{ij}\bold 1)$.

We claim that the scalars $w_{ij}$ are real numbers, so that
$(w_{ij})$ is a real orthogonal $n\times n$ matrix.  Indeed,
if $J_n$ is the standard involution of $A^n$ defined in
(7.13.4) then we have
$$
\align
J_n\theta(\omega) &= \theta(\omega^\#),\qquad\text{and}\\
J_n\theta^\prime(\omega) &= \theta\prime(\omega^\#).
\endalign
$$
It follows that
$$
WJ_n\theta(\omega) = W\theta(\omega^\#) = \theta^\prime(\omega^\#)
= J_n\theta^\prime(\omega) = J_nW\theta(\omega),\qquad \omega\in \Omega^1,
$$
and hence $J_nW = WJ_n$.  Noting that the action of $J_n$ on
scalar operator matrices obeys
$$
J_n(\lambda_{ij}\bold 1)J_n^{-1} = (\overline{\lambda}_{ij}\bold 1),
$$
we conclude that $\overline{w}_{ij} = w_{ij}$ for all $i,j$, as asserted.

Finally, since $\theta^\prime(dx)=W\theta(dx)$, we have
$$
[x_k^\prime,x] = \sum_{j=1}^nw_{kj}[x_j,x] = [\sum_j w_{kj}x_j,x]
$$
for all $x\in A$ it follows that each operator
$$
x_k^\prime - \sum_{k=1}^n w_{kj}x_j
$$
must be a scalar multiple of $\bold 1$.  Since
$\rho(x_k) = \rho(x_k^\prime) = 0$ and $\rho$ is a state, the
scalars must be $0$ and we have the asserted relation
$$
x_k^\prime = \sum_{k=1}^b w_{kj}x_j, \qquad 1\leq k\leq n.
$$
That completes the proof of Theorem 7.6\qed
\enddemo

Theorem 7.6 is applied to the classification of elliptic
operators as follows.

\proclaim{Theorem 7.26}
Let $L\in \Cal D(A,\rho)$ be an exact elliptic operator.  Then there
is a unique momentum space $(P,\<\cdot,\cdot\>)$ and a unique
skew-adjoint operator $v$ such that $\delta(v) = v$, $\rho(v)=0$, and
$$
L = \Delta + [v,\cdot], \tag{7.26.1}
$$
$\Delta$ denoting the Laplacian of $(P,\<\cdot,\cdot\>)$.

Conversely, any operator $L$ of the form (7.26.1) is an exact elliptic
operator.
\endproclaim
\demo{proof}
Since $L$ is an exact elliptic operator its symbol $\sigma_L$ satisfies
condition (iii) of Theorem 6.9, hence
$$
2g = -\sigma_L
$$
defines a KMS metric on $\Omega^2$.  By Theorem 7.6A there is a unique
momentum space $(P,\<\cdot,\cdot\>)$ whose Laplacian $\Delta$ satisfies
$$
2g = -\sigma_\Delta.
$$
The symbol of $L - \Delta$ therefore vanishes, so that $D = L-\Delta$
is a symmetric derivation for which $\rho\circ D = 0$.  We can find
a skew-adjoint operator $v\in A$ for which $D(x) = [v,x]$, and the
condition $\rho\circ D = 0$ implies that $v$ commutes with the density
matrix of $\rho$.  Replacing $v$ with $v-\rho(v)\bold 1$, we obtain
(7.26.1).

The converse was established in Corollary 6.13.
\enddemo

Exact elliptic operators are classified in terms of their momenta
and potentials as follows.

\proclaim{Theorem 7.27}
Let $L_1$, $L_2$ be two exact elliptic operators
in $\Cal D(A,\rho)$ with momentum spaces
$(P_k,\<\cdot,\cdot\>_k)$ and natural
decompositions $L_k = \Delta_k + [v_k,\cdot]$ , $k=1,2$.  Then
$L_1$ and $L_2$ are conjugate iff there is a $*$-automorphism
$\alpha$ of $A$ such that $\rho\circ\alpha = \rho$, and
$$
\align
\alpha(P_1) &= P_2,\tag{7.27.1}\\
\<\alpha(p),\alpha(q)\>_2 &= \<p,q\>_1,\qquad p,q\in P_1\tag{7.27.2}\\
\alpha(v_1) &= v_2.  \tag{7.27.3}
\endalign
$$
\endproclaim
\demo{proof}
Assume that $L_1$ and $L_2$ are conjugate.  Thus there is a
$*$-automorphism $\alpha$ of $A$ such that $\rho\circ\alpha = \rho$
and $\alpha\circ L_1 = L_2\circ \alpha$.  Since $\alpha$ preserves $\rho$,
it induces a unitary operator on $L^2(A,\rho)$ in the natural way:
$$
U_\alpha:x\mapsto \alpha(x),\qquad x\in A,
$$
and we have $U_\alpha L_1 U_\alpha^{-1} = L_2$.  Hence
$U_\alpha(L_1 + L_1^*)U_\alpha^{-1} = L_2 + L_2^*$.  Since
$\Delta_k = L_k + L_k^*$ it follows that $\alpha\Delta_1\alpha^{-1}
= \Delta_2$.  Hence $\alpha$ must induce an isomorphism of momentum
spaces in the sense specified in (7.27.1) and (7.27.2).

Setting $D_k = [v_k,\cdot]$, we have
$$
\alpha\circ D_1\circ\alpha^{-1} = \alpha\circ(L_1-\Delta_1)\circ\alpha^{-1}
= L_2 - \Delta_2 = D_2,
$$
so that for all $x$ in $A$ we have $\alpha([v_i,\alpha^{-1}(x)] = [v_2,x]$.
The latter implies that $\alpha(v_1) - v_2$ is a scalar; and since
$\rho(v_2) = \rho(\alpha(v_1)) = \rho(v_1) = 0$, we deduce (7.27.3).

The proof of the converse is similar\qed
\enddemo

\remark{Remarks}
Of course, the $*$ automorphism $\alpha$ can be implemented by
a unitary operator.

Notice too that the three invariants in
(7.26.1)--(7.26.3) that classify exact elliptic operators can be
specified independently of one another.  For example, starting
with any such operator $L$ one may scale the metric $\<\cdot,\cdot\>$
on $P$ to obtain a family of mutually non-conjugate operators, or
one can scale the potential $v$ to obtain a second family of
mutually non-conjugate operators.

It also makes sense to speak of elliptic operators that are ``free"
in the sense that their potential operator $v$ is zero.
\endremark

\subheading{8.  Applications}
The preceding results lead directly to a new classification
of \esg s.  For that one considers pairs $(\alpha,\omega)$ consisting
of an \esg\ $\alpha = \{\alpha_t: t\geq 0\}$ acting on
$\Cal B(H_\alpha)$ and a normal
state $\omega$ of $\Cal B(H_\alpha)$ satisfying
$$
\omega\circ\alpha_t = \omega,\qquad t\geq 0.
$$
Two such pairs $(\alpha,\omega)$ and $(\alpha^\prime,\omega^\prime)$
are said to be {\it conjugate} if there is a $*$-isomorphism
$\theta:\Cal B(H_\alpha) \to \Cal B(H_{\alpha^\prime})$ such that
$$
\align
\theta\circ \alpha_t &= \alpha^\prime\circ \theta,
\qquad \text{and}\tag{8.1.1}\\
\omega^\prime\circ\theta &= \omega.  \tag{8.1.2}
\endalign
$$

\proclaim{Proposition 8.2}
Assume that $\theta$ is a $*$-isomorphism satisfying
 (8.1.1) and (8.1.2) and let
$p_0$ (resp. $p_0^\prime$) be the support projection of
$\omega$ (resp. $\omega^\prime$).  Let $H_0=p_0H_\alpha$,
$H_0^\prime = p_0^\prime H_{\alpha^\prime}$ and let
$(\phi,\rho)$, $(\phi^\prime,\rho^\prime)$ be the
Markov semigroups obtained by compressing $\alpha$,
$\alpha^\prime$ to $\Cal B(H_0)$, $\Cal B(H_0^\prime)$
as in \S 1:
$$
\align
\phi_t(a) &= p_0\alpha_t(ap_0)p_0,\\
\phi^\prime_t(b) &=p_o^\prime\alpha^\prime_t(bp_0)p_0^\prime,
\qquad t\geq 0,
\endalign
$$
where $\rho$, $\rho^\prime$ are obtained from $\omega$,
$\omega^\prime$ by restriction.  Then the restriction
$\theta_0$ of $\theta$ to $\Cal B(H_0)$ implements a conjugacy
of $(\phi,\rho)$ and $(\phi^\prime,\rho^\prime)$.
\endproclaim
\demo{proof}
The argument is straightforward once one observes that
the condition $\omega^\prime\circ\theta = \omega$ implies that
$\theta(p_0) = p_0^\prime$ and $\rho^\prime\circ\theta_0=\rho$\qed
\enddemo

\proclaim{Definition 8.3}
An \esg\ pair $(\alpha,\omega)$ is said to be of finite
type if the support projection of $\omega$ is finite
dimensional.
\endproclaim

If $(\alpha,\omega)$ is of finite type and $(\phi,\rho)$
is its natural Markov semigroup obtained by compression
as in \S 1
then we may consider the generator $L$ of $\phi$,
$$
\phi_t = \exp(tL),\qquad t\geq 0.  \tag{8.4}
$$
$L$ is an elliptic operator in $\Cal D(\Cal B(H_0),\rho)$ and
it is sensible to ask if $L$ is exact.

\proclaim{Definition 8.5}
$(\alpha,\omega)$ is called exact if it is of finite type
and the generator $L$ of its associated Markov
semigroup is an exact elliptic operator.
\endproclaim

Notice that if $(\alpha,\omega)$ is exact
then the generator $L$ of (8.4) decomposes
uniquely as in Theorem 7.6,
$$
L = \Delta + [v,\cdot]
$$
where $\Delta$ is the Laplacian of a unique momentum space
$(P,\<\cdot,\cdot\>)$ and $v$ is a skew-adjoint operator
satisfying $\rho(v)=0$ and $\rho\circ[v,\cdot] = 0$.  The
triple $(P,\<\cdot,\cdot\>,v)$ is called the {\it dynamical
invariant} of $(\alpha,\omega)$.  It is truly an invariant
because of the following result.

\proclaim{Theorem 8.6}
If two exact \esg\ pairs $(\alpha,\omega)$ and
$(\alpha^\prime,\omega^\prime)$ are conjugate
then their dynamical invariants $(P,\<\cdot,\cdot\>,v)$ and
$(P^\prime,\<\cdot,\cdot\>^\prime,v^\prime)$ are isomorphic in
the sense that there is a $*$-isomorphism
$\alpha:\Cal B(H_0)\to\Cal B(H_0^\prime)$ which satisfies the
conditions of Theorem 7.27:
$$
\align
\omega^\prime\circ\alpha &= \omega  \\
\alpha(P) &= P^\prime\\
\<\alpha(p),\alpha(q)\>^\prime &= \<p,q\>,\qquad p,q\in P\\
\alpha(v) &= v^\prime.
\endalign
$$
\endproclaim

\demo{proof}
Given the preceding remarks, the argument is a simple variation
of what was done in the proof of Theorem 7.27\qed
\enddemo

\remark{Remarks}
We have already pointed out in \S 1 that the converse of Theorem
8.6 is false in general, but that it is true if both $\alpha$
and $\alpha^\prime$ are {\it minimal}. More generally, we can
paraphrase a recent dilation theorem of B. V. R. Bhat \cite{6}
in our context as follows.
\endremark

\proclaim{Theorem 8.7} Let $(\phi,\rho)$ be a Markov semigroup
acting on a separable Hilbert space $H_0$.  Then there is a
minimal \esg\ $\alpha$ acting on a Hilbert space $H\supseteq H_0$
and a normal state $\omega$ on $\Cal B(H)$ such that
$\omega\circ\alpha_t = \omega$ for every $t\geq 0$ and such that
$(\phi,\rho)$ is obtained from $(\alpha,\omega)$ by compression
as described above.

If $(\phi^\prime,\rho^\prime)$ is another Markov semigroup which
is conjugate to $(\phi,\rho)$ and which gives rise to a minimal
\esg\ pair $(\alpha^\prime,\omega^\prime)$ then $(\alpha,\omega)$
and $(\alpha^\prime,\omega^\prime)$ are conjugate.
\endproclaim

\remark{Remark 8.8} Taken together, Theorems 8.6 and 8.7 imply that
exact {\it minimal} \esg\ pairs $(\alpha,\omega)$ are completely classified
up to conjugacy by their dynamical invariants $(P,\<\cdot,\cdot\>, v)$.
More explicitly, suppose we start with a pair $(A,\rho)$ consisting
of a faithful state $\rho$ on $A=\Cal B(H_0)$, where $H_0$ is a
finite dimensional Hilbert space.  Choose an arbitrary dynamical
invariant $(P,\<\cdot,\cdot\>,v)$ in $(A,\rho)$ and form its
corresponding elliptic operator $L\in \Cal D(A,\rho)$:
$$
L(x) = \Delta(x) + [v,x],\qquad x\in A,
$$
$\Delta$ being the Laplacian of $(P,\<\cdot,\cdot\>)$.  One
exponentiates $L$ to obtain $\phi_t = \exp(tL)$ and
a Markov semigroup $(\phi,\rho)$.  Finally, one obtains
a minimal \esg\ pair $(\alpha,\omega)$ by the dilation procedure
of Theorem 8.7.
Such an $(\alpha,\omega)$ is an exact pair whose dynamical
invariant is isomorphic to $(P,\<\cdot,\cdot\>,v)$.  Every exact
pair arises in this way, and two exact pairs are conjugate
if and only if their corresponding dynamical invariants are
isomorphic.

We infer from these remarks that the dynamical invariant
$(P,\<\cdot,\cdot\>,v)$ consists of infinitesimal structures that
completely determine the associated flow $(\alpha,\omega)$.  Thus
dynamical invariants are analogous to the differential equations that
govern the behavior of classical mechanical systems.
There are obvious implications for the classification of type $I$
primary histories.

Finally, it is possible to determine which triples $(P,\<\cdot,\cdot\>,v)$
lead to ``mixing" properties for their associated \esg s.  By that we
mean that all limits of the form
$$
\lim_{t\to\infty}\phi(\alpha_t(a)),\qquad a\in \Cal B(H_\alpha),
\phi\in \Cal B(H_\alpha)_*
$$
should exist.  These issues will be taken up elsewhere.
\endremark

\vfill
\pagebreak

\end